\documentclass[3p,preprint]{elsarticle}

\usepackage{amsmath}
\usepackage{amsfonts}
\usepackage{amssymb,latexsym}
\usepackage{amsthm}
\usepackage{graphicx}
\usepackage{graphics}
\usepackage{epsfig}
\usepackage{multirow}
\usepackage{graphicx}
\usepackage{adjustbox}

\usepackage{algorithm}
\usepackage{algpseudocode}

\usepackage{url}

\usepackage{color}
\usepackage{newlfont}
\usepackage{subfigure}
\usepackage{verbatim}
\usepackage{rotating}
\usepackage{multirow}
\usepackage{units}
\usepackage{bm}
\usepackage[english]{babel}
\usepackage[utf8]{inputenc}
\usepackage{algorithm}
\usepackage{booktabs}
\usepackage{caption}
\usepackage{hyperref}

\usepackage{textcomp}
\DeclareGraphicsExtensions{.jpg,.jpeg,.pdf,.png,.mps,.eps,.ps,.bmp}
\usepackage{eurosym}
\usepackage[toc,page]{appendix}
\usepackage{graphicx}
\usepackage{subfigure}
\usepackage{graphicx,epstopdf}

\usepackage{caption}
\usepackage{verbdef}
\usepackage{listings}
\usepackage{fancyvrb}
\usepackage{float}
\usepackage{pdfpages}
\usepackage{verbatim}
\usepackage{lmodern}
\usepackage{textcomp}
\usepackage{titlesec}

\DeclareMathOperator*{\argmax}{argmax}

\newcommand{\mK}{\mathsf{K}}

\journal{Journal of Computational and Applied Mathematics}

\begin{document}

\begin{frontmatter}



\title{Bayesian Approach for Radial Kernel Parameter Tuning}




\author[address-TO,address-GNCS]{Roberto Cavoretto\corref{corrauthor}} 
\ead{roberto.cavoretto@unito.it}

\author[address-TO,address-GNCS]{Alessandra De Rossi}
\ead{alessandra.derossi@unito.it}

\author[address-TO,address-GNCS]{Sandro Lancellotti}
\ead{sandro.lancellotti@unito.it}

\cortext[corrauthor]{Corresponding author}

\address[address-TO]{Department of Mathematics \lq\lq Giuseppe Peano\rq\rq, University of Torino, via Carlo Alberto 10, 10123 Torino, Italy}

\address[address-GNCS]{Member of the INdAM Research group GNCS}

\begin{abstract}
In this paper we present a new fast and accurate method for Radial Basis Function (RBF) approximation, including interpolation as a special case, which enables us to effectively find  the optimal value of the RBF shape parameter. In particular, we propose a statistical technique, called Bayesian optimization, that consists in modelling the error function with a Gaussian process, by which, through an iterative process, the optimal shape parameter is selected. The process is step by step self-updated resulting in a relevant decrease in search time with respect to the classical leave one out cross validation technique.
Numerical results deriving from some test examples and an application to real data show the performance of the proposed method.
\end{abstract}




\end{frontmatter}



\section{Introduction}
Kernel-based techniques are effective procedures commonly used to interpolate and approximate scattered data sets. They owe their popularity in the field of approximation theory, notably in meshfree approximation, not only to the capability to be used in high dimensions but also to the attainable accuracy of the approximants. Such methods are obtained by a weighted sum of some radial kernels or Radial Basis Functions (RBFs) \cite{fas07}, which depend on the so-called {\sl shape parameter}. The latter significantly affects the accuracy of the radial kernel method, and so its detection is paramount. In particular, there exists a trade-off between accuracy and stability (see \cite{wen05}). This issue drove authors to choose the parameter by a trial and error approach ending up proposing empirical results with a whopping waste of time and computational resources or some ad-hoc criteria \cite{Franke, Hardy} that, sometimes, could be non-optimal. Other systematic approaches for the best value or an optimal interval have been proposed in \cite{Biazar,Scheuerer} with the benefit of being generally applicable but computationally expensive. Notice that the issue of selecting an optimal scale parameter is a current research topic in several fields of applied mathematics and scientific computing (see e.g. \cite{cav21a,cav21,lin22}). 

The purpose of this work is to propose a flexible and suitable technique applicable to approximation and interpolation capable to head to a good error estimate without spending much time in trials. The key idea is to apply a well-known statistical technique, called {\sl Bayesian Optimization} (BO) \cite{Practical}, to search for the optimal shape parameter. This technique, often considered in machine learning  for optimization of black-box or difficult-to-evaluate functions, can be used in problems of hyper-parameter tuning to avoid the computation and evaluation of the approximant for those parameters which are far from being optimal. Accordingly, this leads to a remarkable saving of computational time during the exploration phase in the parameter domain.

To complete this study, we compare BO results with those obtained by using  Rippa's version of the Leave One Out Cross Validation (LOOCV) scheme \cite{rip99} that, given $n$ points, reduces the computational cost from $n^4$ to $n^3$ with respect of the standard implementation previously introduced in \cite{Golberg}. It is a popular and well-established strategy  developed for selecting a good RBF shape parameter for data interpolation and later popularised in the field of applied sciences for training performance evaluation \cite{Trahan, Uddin}. 


The rest of paper is organized as follows. In Section \ref{rbf}, RBF interpolation and least squares RBF approximation are briefly stated. In Section \ref{optimizers}, LOOCV and BO are described. In Section \ref{algorithms} the proposed algorithms are explained in detail. In Section \ref{exp} numerical results on some test examples are shown, while in Section \ref{real_data_example}  an application to a real dataset is considered. Finally, Section \ref{conclusion} concludes the paper. 

\section{RBF Theory} \label{rbf} 

In the present section we introduce the problem, the basic theory about RBF interpolation and approximation, and the notation we will use throughout the paper. Moreover, we will also discuss the reason that inspired this work, i.e., the search of the shape parameter associated with the RBF.

 \subsection{RBF Interpolation} \label{rbf_interp}
 
Suppose to have a set of distinct data points or data locations $ X = \{  \boldsymbol{x}_i, i = 1,  \ldots , n \}$ arbitrarily distributed on a domain $ \Omega \subseteq \mathbb{R}^{d}$, with an associated set $ F = \{ f_i = f(\boldsymbol{x}_i) ,i=1, \ldots, n \}$ of data values, which are obtained by sampling some (unknown) function  $f: \Omega \longrightarrow \mathbb{R}$ at the nodes $ \boldsymbol{x}_i$, the scattered data interpolation problem consists in finding a function $P_f: \Omega \longrightarrow \mathbb{R}$ such that it matches the measurements at the corresponding locations, i.e. 
$$
P_f\left( \boldsymbol{x}_i\right)=f_i, \quad i=1, \ldots, n.
$$

We now suppose to have a univariate function $ \varphi: [0, \infty) \to \mathbb{R}$, known as RBF, which depends on a shape parameter $\varepsilon > 0$ providing, for $\boldsymbol{x},\boldsymbol{z}
\in \Omega$, the real symmetric strictly positive definite kernel (see e.g. \cite{wen05})
\begin{equation*} 
\kappa_\varepsilon(\boldsymbol{x},\boldsymbol{z}) = \varphi(\varepsilon ||\boldsymbol{x}-\boldsymbol{z}||_2): = \varphi(\varepsilon r).
\end{equation*}

The kernel-based interpolant $P_f$ can be written as
\begin{equation} \label{ker_interp}
P_f\left( \boldsymbol{x}\right)= \sum_{k=1}^{n} c_k \kappa_\varepsilon \left( \boldsymbol{x} , \boldsymbol{x}_k  \right), \quad \boldsymbol{x} \in \Omega,
\end{equation}
whose coefficients $c_k$ are the solution of the linear system
\begin{equation} \label{linsys}
\mK_\varepsilon \boldsymbol{c} = \boldsymbol{f},
\end{equation}
where $  \boldsymbol{c}= \left(c_1, \ldots,
c_n\right)^{\intercal}$, $  \boldsymbol{f} =\left(f_1, \ldots , f_n\right)^{\intercal}$, and $(\mK_{\varepsilon})_{ik}= \kappa_\varepsilon \left( \boldsymbol{x}_i , \boldsymbol{x}_k  \right)$, $i,k=1, \ldots, n$. Since $\kappa_\varepsilon$ is a symmetric and strictly positive definite kernel, the system \eqref{linsys} has exactly one solution \cite{Fasshauer15}. Additionally, for the kernel $\kappa_\varepsilon$ there exists the so-called \emph{native space}, which is a Hilbert space ${\cal N}_{\kappa_\varepsilon}(\Omega)$ with inner product $(\cdot,\cdot)_{{\cal N}_{\kappa_\varepsilon}(\Omega)}$ in which the kernel $\kappa_\varepsilon$ is reproducing, i.e., for any $f \in {\cal N}_{\kappa_\varepsilon}(\Omega)$ we have the identity $f(\boldsymbol{x}) = (f,\kappa_\varepsilon(\cdot,\boldsymbol{x}))_{{\cal N}_{\kappa_\varepsilon}(\Omega)}$, with $\boldsymbol{x}\in \Omega$. Then, if we introduce a pre-Hilbert space $H_{\kappa_\varepsilon}(\Omega)= \mbox{span}\{\kappa_\varepsilon(\cdot,\boldsymbol{x}),$ $\boldsymbol{x} \in \Omega\}$, with reproducing kernel $\kappa_\varepsilon$ and equipped with the bilinear form $(\cdot,\cdot)_{H_{\kappa_\varepsilon}(\Omega)}$, the native space ${\cal N}_{\kappa_\varepsilon}(\Omega)$ of $\kappa_\varepsilon$ is its completion with respect to the norm $||\cdot||_{H_{\kappa_\varepsilon}(\Omega)}=\sqrt{(\cdot,\cdot)_{H_{\kappa_\varepsilon}(\Omega)}}$. In particular, for all $f \in {H_{\kappa_\varepsilon}(\Omega)}$ we have $||f||_{{\cal N}_{\kappa_\varepsilon}(\Omega)}=||f||_{H_{\kappa_\varepsilon}(\Omega)}$. Error bounds for kernel interpolants \eqref{ker_interp}, expressed in terms of the well-known \textsl{power function} $P_{\kappa_\varepsilon, X}$, can be found, for instance, in \cite{fas07}. For more refined error estimates, see e.g. \cite{wen05}. For the sake of completeness, however, we also observe that not all radial kernels are depending on a shape parameter. Indeed, there exist shape parameter free kernels like polyharmonic splines that are very popular and commonly used tools in the RBF community \cite{mir21}.

 \subsection{Least Squares RBF Approximation} \label{rbf_approx}

So far we have looked only at interpolation. Nevertheless, sometimes it makes more sense to approximate the given data by a least squares RBF approximant. Indeed, if the data are subjected to noise, or there are so many data locations that efficiency motivations force us to approximate from a space spanned by fewer basis functions than data locations.

In this subsection, we consider a more general setting where we still sample the given function $f$ at the data point set $X$ but here we also introduce a set of centers $\tilde{X}=\{\tilde{\boldsymbol{x}}_i, i=1,\ldots,m\}$ at which the kernels are centred. Generally, we have $m \leq n$ and the case $m = n$ with $\tilde{X} = X$ recovers the RBF interpolation setting discussed in Subsection \ref{rbf_interp}.

So the RBF approximant may be expressed as
\begin{align}\label{ker_approx}
\tilde{P}_f\left( \boldsymbol{x}\right)= \sum_{k=1}^{m} \tilde{c}_k \kappa_\varepsilon \left( \boldsymbol{x} , \tilde{\boldsymbol{x}}_k  \right), \quad \boldsymbol{x} \in \Omega,
\end{align}
the coefficients $\tilde{c}_k$ being determined as the least squares solution of the linear system
\begin{equation} \label{linsys_hat}
\tilde{\mK}_{\varepsilon} \tilde{\boldsymbol{c}} = \boldsymbol{f},
\end{equation}
which is obtained by minimizing $||\tilde{P}_f - f||_2^2$, where $(\tilde{\mK}_{\varepsilon})_{ik}= \kappa_\varepsilon \left( \boldsymbol{x}_i , \tilde{\boldsymbol{x}}_k  \right)$, $i=1, \ldots, n$, $k=1, \ldots, m$, and $\boldsymbol{f} =\left(f_1, \ldots , f_n\right)^{\intercal}$. This approximation problem has a unique solution if the $n \times m$ collocation matrix $\tilde{\mK}_\varepsilon$ has full rank. Now, if the centers $\tilde{X}$ are a subset of the data locations $X$, then the matrix $\tilde{\mK}_\varepsilon$ has full rank provided that the kernel is symmetric and strictly positive definite \cite{fas07}.

As the accuracy of the fit strongly depends from the choice of the shape parameter $\varepsilon$, see e.g. \cite{cav22,Fornberg-Wright04,Larsson-Fornberg05,Marchetti}, in the following section we describe two different techniques for the search of the shape parameter. The former is a classical approach, known as the LOOCV, while the latter is the BO which represents our proposal to reduce the computational expense and avoids some drawbacks related to LOOCV. We would like also to highlight that during the experiments we used  Rippa's version of the LOOCV that can only be applied  on interpolation problems. Indeed, although the standard LOOCV technique could be utilized in approximation settings, we did not perform experiments due to the high theoretical computational cost.

\section{Optimizers} \label{optimizers}
In the field of kernel-based approximation, the search of the shape parameter is one of the topics that most attracts the attention of researchers due to the strong dependence of the fit accuracy on it. In this section we briefly describe a well-known technique, the LOOCV, and a statistical procedure, the BO, as an alternative approach to avoid the LOOCV deficiencies. 

\subsection{Leave One Out Cross Validation}

A popular strategy for estimating the RBF shape parameter $\varepsilon$ based on the given data set $(X,F)$ is the LOOCV method. In this technique an optimal value of $\varepsilon$ is selected by minimizing a cost function that collects the errors for a sequence of partial fits to the data. To estimate the unknown true error, we split the data into two parts: a \emph{training} data set consisting of all data, except for one, to obtain a partial fit, and a \emph{validation} data set that contains the single remaining datum used to compute the error. After repeating in turn this procedure for each of the given data, the result is a vector of error estimates and the cost function is used to determine the optimal value of $\varepsilon$, see \cite{fas07b}.

The LOOCV is a technique for the search of the optimal value of the RBF shape parameter $\varepsilon$. It consists in evaluating, for each $\varepsilon$ and for each  $ j \in \lbrace 1, \dots, n \rbrace$, the error 
$$e_j(\varepsilon) = f(\boldsymbol{x}_j) - P_f^{j}( \boldsymbol{x}_j)$$ 
at the validation point $\boldsymbol{x}_j$ that is not used to construct the partial RBF interpolant
\begin{align} \label{ker_pinterp}
	P_f^j( \boldsymbol{x}) = \sum_{k=1,\ k \neq j}^{n} c_k \kappa_\varepsilon(\boldsymbol{x},\boldsymbol{x}_k).
\end{align} 
The latter is fitted on the training data point set $X_j = X \setminus \lbrace x_j \rbrace$ and the data values  $F_j = F \setminus \lbrace f_j \rbrace$, while the coefficients $c_k$ in \eqref{ker_pinterp} are determined by interpolating only the set $X_j$, i.e.,  
\begin{align*}
P_f^j( \boldsymbol{x}_k) = f(\boldsymbol{x}_k), \quad k=1,\ldots,j-1,j+1,\ldots,n.
\end{align*}
The optimal value of $\varepsilon$ is found as
\begin{align*} 
\varepsilon^* = \mbox{argmin}_{\varepsilon} ||\boldsymbol{e}(\varepsilon)||, \qquad \boldsymbol{e} = (e_1,\ldots,e_{n})^{\intercal},
\end{align*}
where $||\cdot||$ is any norm used in the minimization problem, for instance, the $\infty$-norm.

Since this LOOCV implementation is quite expensive, the error computation can be simplified by using the rule proposed by Rippa in \cite{rip99}
\begin{equation*}
    e_j(\varepsilon) = \frac{c_j}{(\mK_\varepsilon^{-1})_{jj}}.
\end{equation*}
where $c_j$ is the $j$th coefficient of the solution vector $\boldsymbol{c} = \mK_\varepsilon^{-1}\boldsymbol{f}$ in \eqref{linsys}, and $(\mK_\varepsilon^{-1})_{jj}$ is the $j$th diagonal element of the inverse of the \emph{full} RBF matrix $\mK_\varepsilon$. Notice that this formulation needs to only solve a single linear system, considering the entire data set $X$ and thus, avoiding the solution of $n$ interpolation problems on $n-1$ points, allows to reduce the computational cost from $n^4$ to $n^3$. It is worth noting that when least square approximation is performed, the matrix $\mK_\varepsilon$ is non-invertible and therefore Rippa's formulation is not applicable.

It follows immediately that the optimal value $\varepsilon^*$ for the shape parameter is the one that minimizes the error function $Er(\varepsilon)$ defined as follows:
\begin{equation}\label{error_function}
Er(\varepsilon) = \max_{j=1, \dots, n} \Bigg|   \frac{c_j}{(\mK_\varepsilon^{-1})_{jj}} \Bigg|.
\end{equation}

In order to find $\varepsilon^*$ (or, at least, a good approximation of it), a finite set of equally spaced values between $0$ and a large enough $\varepsilon_{max}$ are exploited to evaluate the error function \eqref{error_function}. 
The LOOCV technique is formalized by imposing to evaluate the error function also for that value in the discrete set that does not lead to a good result. 
Moreover, applying this scheme, in general, is not possible to attain the global minimum due to the discrete search. Another drawback related to the evaluation of the error is the inversion of the interpolation matrix which could lead to additional avoidable instabilities during the computation.  
A possible extension of LOOCV could be the addition of a univariate optimiser to direct the search of the  $\varepsilon$ parameter. Such a solution would certainly bring an improvement in computation time performance but it would also be subject to the choice of the initial parameter value from which to start the search. Furthermore, in order to apply the LOOCV with an optimiser to the approximation case, one would have to dispense with the use of Rippa's formula. This fact would lead to a high computational cost and would not justify its use.
The aim of our work is to propose an alternative technique that saves computational time, avoids the evaluation of bad values of $\varepsilon$ and is able to conduct a continuous search in the parameter space, thus producing a better approximation of $\varepsilon^*$ and avoiding all related possible failures. The reader can find a comprehensive treatment with examples of the drawbacks related to the search of the shape parameter in \cite[Chapter 17]{fas07}.

\subsection{Bayesian Optimization}
\label{bo}
When it comes to find a global maximizer for an unknown or difficult-to-evaluate function $g$ on some bounded set $X$, the Bayesian optimization \cite{Mockus_1978} is an elegant approach to carry out the search. Very popular in machine learning, the Bayesian optimization is an iterative technique based on a simple principle: exploiting all the available resources. It consists in building a probabilistic model of $g$, called  {\sl surrogate model}, and using it to help directly the sampling point in $X$, by means of an acquisition function, where the target function will be evaluated. As an iteration is made, the distribution is first updated and then used in the next iteration. Even though there is a computation for the selection of the next point to evaluate, when evaluations of $g$ are expensive, it turns out that the computation of a better point is motivated by reaching the maximum in a few iterations, as in the case of the error function of some expensive training machine learning algorithms like multi-layer neural networks (see \cite{review,Bishop}).


In this section, we briefly review the Bayesian optimization technique referring the reader to \cite{Brochu} for a more detailed description.

\vskip 0.1cm
\textbf{Gaussian Processes.} A Gaussian Process (GP) is a collection of random variables such that any subsets of them have a joint Gaussian  distribution. It is completely specified by a mean function $m : \mathcal{X} \rightarrow \mathbb{R}$ and a positive definite covariance function $k : \mathcal{X} \times \mathcal{X} \rightarrow \mathbb{R}$ where $\mathcal{X} \subseteq \mathbb{R}$ (see \cite{Rasmussen} for further details).

Gaussian processes are the most common choice for the surrogate model for Bayesian optimization due to the low evaluation cost and to the ability to incorporate prior beliefs about the objective function. When we model the target function with a Gaussian process as $g(x) \sim \mathcal{GP} \big( m(x), k(x, x') \big)$, we are imposing the following conditions:
\begin{itemize}
    \item $\mathbb{E}\big [ g(x) \big] = m(x)$;
    \item $\mathbb{E}\big [ \big( g(x) - m(x) \big)  \big( g(x') - m(x') \big) \big] = k(x, x')$.
\end{itemize}
When it comes to make a prediction given by some observations, the assumption of joint Gaussianity allows retrieving the prediction using the standard formula for the mean and variance of a conditional normal distribution. Notably, suppose to have $s$ observations $\bm{g} = (g(x_1), \dots, g(x_s))^{\intercal}$ on the points $\textbf{x} = (x_1, \dots, x_s)^{\intercal}$  and a new point $\bar{x}$ where we are interested in having a prediction $\bar{g}$ of $g(\bar{x})$. The previous observations $\bm{g}$ and the predicted value $\bar{g}$ are jointly normally distributed:
$$
Pr\begin{pmatrix} \begin{bmatrix} \bm{g} \\ \bar{g} \end{bmatrix} \end{pmatrix} = \mathcal{N} \begin{bmatrix}
\begin{bmatrix} \mu(\textbf{x}) \\ \mu(\bar{x}) \end{bmatrix}, 
\begin{bmatrix} 
K(X,X) \ \ K(\textbf{x}, \bar{x}) \\
K(\textbf{x}, \bar{x})^{\intercal} \ \ k(\bar{x}, \bar{x}) \\
\end{bmatrix}\end{bmatrix},
$$
where $X$ is the $s \times s$ matrix with $(i, j)$-element $(x_i, x_j)$,  $K(X,X)$ is the $s \times s$ matrix with $(i, j)$-element $k(x_i, x_j)$, and
$K(\textbf{x}, \bar{x})$ is a $s \times 1$ vector whose element $i$ is given by $k(x_i, \bar{x})$. Since $Pr(\bar{g} |  \bm{g})$ must also be normal, it follows that:
$$
Pr(\bar{g} | \bm{g}) = \mathcal{N}  \begin{bmatrix} \mu(\bar{x})  K(X, \bar{x})^{\intercal} K(X,X)^{-1} (\bm{g} - \mu(\textbf{x})) , k(\bar{x}, \bar{x}) - K(X, \bar{x})^{\intercal}  K(X, X)^{-1} K(X, \bar{x})
\end{bmatrix}.
$$

\begin{figure}
    \centering
    \includegraphics[width=.33\textwidth]{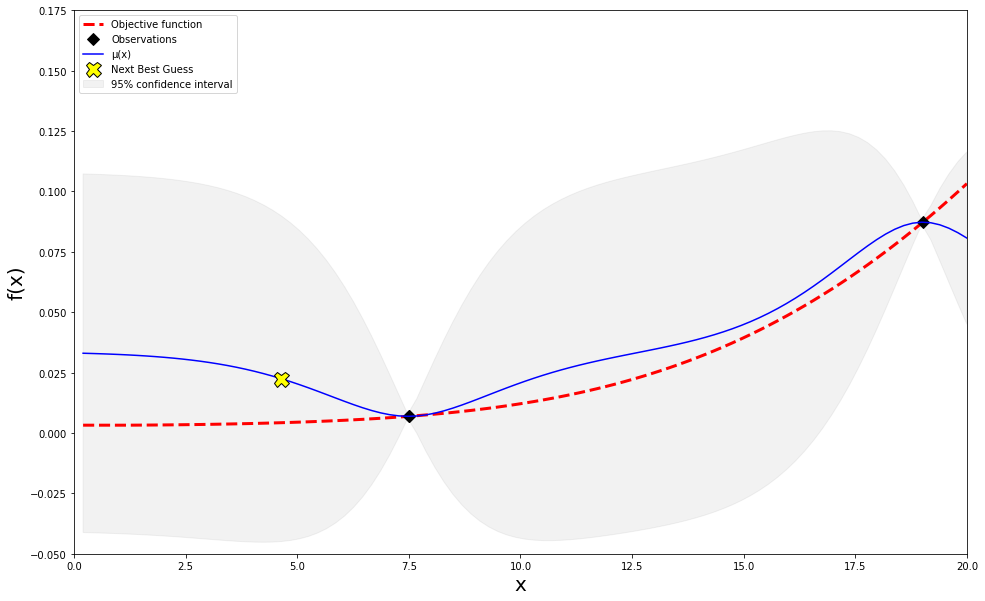}\hfill
\includegraphics[width=.33\textwidth]{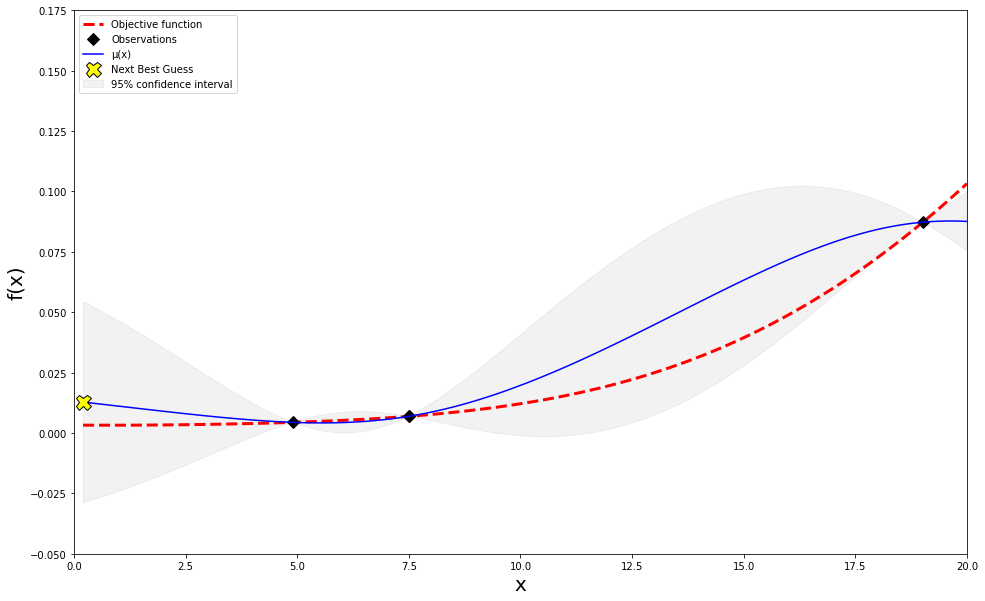}\hfill
\includegraphics[width=.33\textwidth]{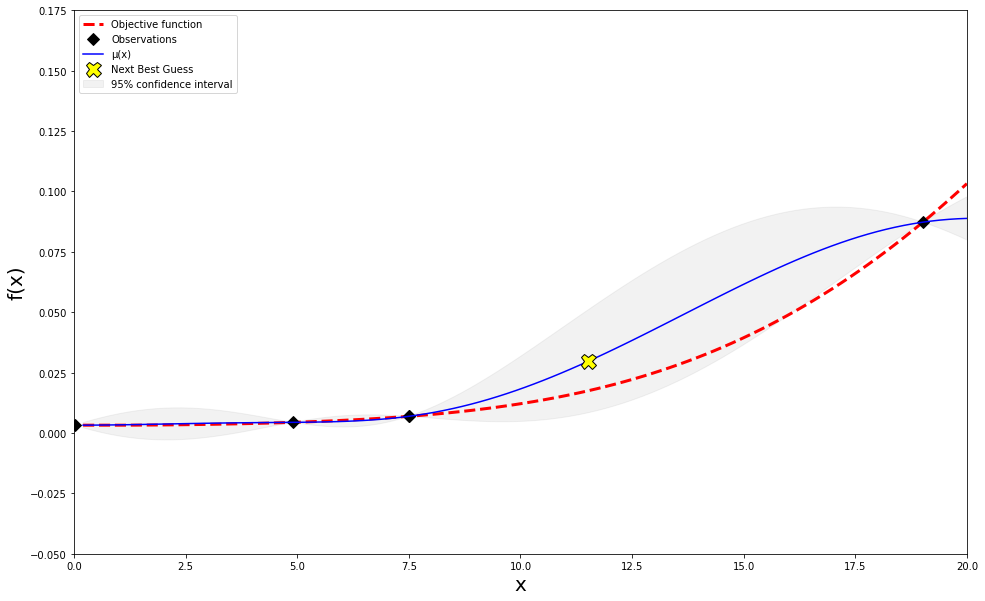}
    \caption{An example of 2 Bayesian optimization steps with 2 starting points. The left image is the snapshot before the first Bayesian iteration while the others are the two following Bayesian steps. The yellow cross is the point selected by the Expected Improvement acquisition function and it is used to update the Gaussian model. As can be seen, at each iteration the variance of the process became smaller and the confidence interval narrows around the objective function.}
    \label{fig:bayesian_steps}
\end{figure}

In this way, it is possible to estimate the distribution, mean and covariance, at any point in the domain. When data locations and data values retrieved by the evaluation of the target function are fed to the model, they induce a posterior distribution over functions which is used for the next iteration as a prior (see Figure \ref{fig:bayesian_steps}). In particular, if a function is modelled by a GP, when we observe a value, we are observing the random variable associated to the point.

\vskip 0.1cm
\textbf{Acquisition Functions.} An acquisition function $a: \mathcal{X} \rightarrow \mathbb{R}$ is a function used to determine the next point to evaluate by the objective function. The chosen point is the one that maximizes this acquisition function, and its evaluation by the objective function is used to update the surrogate model (see Figure \ref{fig:bayesian_steps}). An acquisition function is defined such that high acquisition corresponds to potentially high values of the objective function. There exists a trade-off between exploration and exploitation in the selection of an acquisition function: exploration means selecting points where the uncertainty is high, that is, far from the already evaluated points; exploitation, on the contrary, means selecting those points  close to those already evaluated by the objective function. The most common acquisition functions are:

\begin{itemize}
    \item \textbf{Probability of Improvement}: maximize the probability of improvement over the best current value;
    \item \textbf{Expected Improvement}: maximize the expected improvement over the current best;
    \item \textbf{GP Upper Confidence Bound}: minimize the cumulative regret\footnote[1]{Regret is a performance metric commonly used in Reinforcement Learning. In a maximization setting of a function $g$ it represents the loss in rewards due to not knowing $g$'s maximum points beforehand. If $x^* = \argmax g(x)$, the regret for a point $x$ is $g(x^*) - g(x)$.} over the course of the optimization.
\end{itemize}
The acquisition function we used in this work is the \lq\lq Expected Improvement\rq\rq\ \cite{EI} that takes into account not only the probability of improvement of the candidate point with respect to the previous maximum, but also the magnitude of this improvement.

\begin{figure}
    \centering
    \includegraphics[scale=0.4]{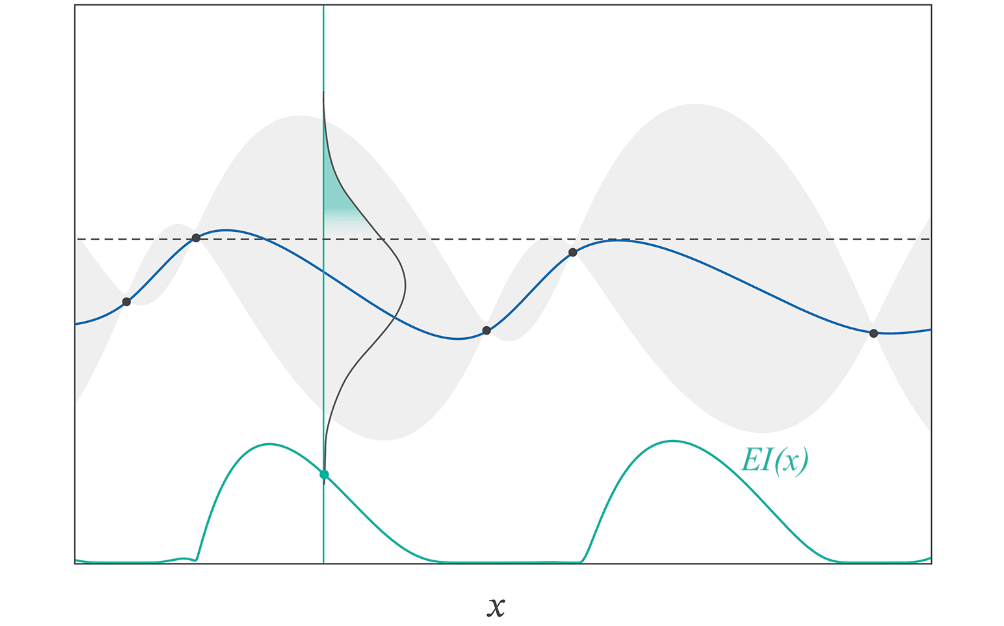}
    \caption{The black points are the evaluation of the objective function while the blue line and the grey shaded area are the mean $\mu(x)$ and the $95\%$ confidence interval of the surrogate model. The teal line represents the Expected Improvement $EI(x)$.}
    \label{fig_EI}
\end{figure}

Suppose that after a number of iterations the current maximum of the objective function is $g(\hat{x})$. Given a new point $x$, the Expected Improvement acquisition function computes the expectation of improvement $g(x) - g(\hat{x})$  over the part of the normal distribution that is above the current maximum (see Figure \ref{fig_EI}):
\begin{equation}
\label{EI_eq}
    EI(x) = \int_{g(\hat{x})}^\infty \big( g^*(x) - g(\hat{x}) \big) 
    \frac{1}{\sqrt{2 \pi} \sigma(x)} e^{-\frac{1}{2}  [(g^*(x) - \mu(x))/\sigma(x)]^2} dg^*(x),
\end{equation}
where $g^*(x)$, $\mu(x)$ and $\sigma(x)$ represent the predicted value by the surrogate model, the expected value and the variance of $x$, respectively.
Solving integral \eqref{EI_eq} leads to the following closed form for the evaluation of the Expected Improvement:
\begin{equation} \label{EI_closed}
    EI(x) = 
    \begin{cases}  
            (\mu(x) - g(\hat{x})) \Phi(Z) + \sigma(x) \phi(Z), & \text{ if } \sigma(x) >0, \\
            0, & \text{ if } \sigma(x) =0, 
    \end{cases}
\end{equation}
where $Z =\frac{\mu(x)-g(\hat{x})}{\sigma(x)}$, while $\phi$ and $\Phi$ are the Probability Density Function and Cumulative Distribution Function of the standard normal distribution $\mathcal{N}(0, 1)$.

An extension of (\ref{EI_closed}) that also trades off exploration and expectation by means of a non-negative parameter $\xi$ was proposed in \cite{Lizotte}: 
\begin{equation} 
\label{EI_closed_trade-off}
    EI(x) = 
    \begin{cases}  
            (\mu(x) - g(\hat{x}) -\xi) \Phi(Z) + \sigma(x) \phi(Z), & \text{ if } \sigma(x) >0, \\
            0, & \text{ if } \sigma(x) =0, 
    \end{cases}
\end{equation}
where $Z =\frac{\mu(x)-g(\hat{x}) -\xi}{\sigma(x)}$. 

\section{Algorithms}
\label{algorithms}
In this section we will explore in detail how Bayesian optimization works when it is applied to an approximation or interpolation process. 

Suppose to have a set of points $X$ for which we know the associated set of data values $F$. 
We consider a subset $\tilde{X} \subseteq X$ of the RBF centers. When the set of centers is strictly contained, we face an approximation problem in the least squares sense; on the other hand, when the equality holds, we fall into the special case of interpolation.

\textbf{Algorithm} \ref{alg_1} is a pseudocode adaptation of the  \em{BayesianOptimization} library \cite{BO}. Specifically, it traces the optimisation process carried out by the \em{optimisation} method of the \em{BayesianOptimisation} class, which in turn uses the methods \em{fit} and \em{predict} of the function \em{GaussianProcessRegressor} of \em{sklearn.gaussian\_process} package \cite{sklearn}. For the details about the  implementation of  \em{GaussianProcessRegressor} we refer the reader to \cite[Algorithm 2.1]{Rasmussen}. 
To apply the Bayesian optimization (see \textbf{Algorithm} \ref{alg_1}) and  evaluate the validation error during the optimization, we first start dividing the sets $X, \tilde{X}$ and $F$ into $X_{train}, X_{val}$,  $\tilde{X}_{train}, \tilde{X}_{val}$ $F_{train}, F_{val}$  in such a way 
that 
\begin{align*}
    |X_{train}| & =  \lfloor0.8 \times |X| \rfloor, \qquad  |X_{val}| = \lceil 0.2\times |X| \rceil, \\ |\tilde{X}_{train}| &= \lfloor 0.8 \times |\tilde{X}| \rfloor, \qquad |\tilde{X}_{val}| = \lceil 0.2\times |\tilde{X}| \rceil, \\ |F_{train}| &= \lfloor 0.8 \times |F| \rfloor, \qquad \  |F_{val}|  = \lceil 0.2\times |F| \rceil.
\end{align*}


\begin{algorithm}
\caption{$\textbf{BO}(X_{train},  \tilde{X}_{train}, F_{train}, X_{val}, g,  \mathcal{X}, a, \xi, nstart,  niter)$}
\hspace*{\algorithmicindent} \textbf{Input:} 
    \vspace{-3mm}
\begin{itemize}
\setlength\itemsep{-0.2em}
    \item[] $X_{train}$: data locations, $\tilde{X}_{train}$: RBF centers, $F_{train}$: data values, $X_{val}$: evaluation points,  $g$: function to maximize, $\mathcal{X}$: parameter search interval, $a$: acquisition function, $\xi$: exploration-exploitation parameter, $nstart$: number of starting iterations, $niter$: number of Bayesian iterations.
\end{itemize}
\begin{algorithmic}

\State  $\bm{\varepsilon} \rightarrow$ sample of $nstart$ random parameter values in $\mathcal{X}$
\State $\textbf{P}_{\textbf{f}_{{\varepsilon}_j}} \rightarrow \textbf{RBF}(X_{train}, \tilde{X}_{train}, F_{train}, X_{val}, {\varepsilon}_j), \qquad j=1, \dots, nstart$  \qquad (call to \textbf{Algorithm} \ref{alg_2})
\State $\bm{g} \rightarrow (g(\textbf{P}_{\textbf{f}_{{\varepsilon}_1}}), \dots, g(\textbf{P}_{\textbf{f}_{{\varepsilon}_{nstart}}}))$

\For{$i = 1 : niter$}
    \State Fit the Gaussian process on  $(\bm{\varepsilon}, \bm{g})$
    \State Evaluate $a$ on a set of random parameter values in $\mathcal{X}$
    \State Select the parameter $\hat{\varepsilon}$ that maximizes $a$
    \State $\textbf{P}_{\textbf{f}_{\hat{\varepsilon}}} \rightarrow \textbf{RBF}(X_{train}, \tilde{X}_{train}, F_{train}, X_{val}, \hat{\varepsilon})$ \qquad (call to \textbf{Algorithm} \ref{alg_2})
    \State $\bm{\varepsilon} \rightarrow \bm{\varepsilon} \cup \hat{\varepsilon} $
    \State $\bm{g} \rightarrow \bm{g} \cup g({\textbf{P}_{\textbf{f}_{\hat{\varepsilon}}}})) $
\EndFor
\State $\varepsilon^* \rightarrow \argmax \bm{g}$
\end{algorithmic}
\hspace*{\algorithmicindent} \textbf{Output:} 
    \vspace{-3mm}
\begin{itemize}
\setlength\itemsep{-0.2em}
    \item[] $\varepsilon^*$: Best shape parameter.
\end{itemize}
\label{alg_1}
\end{algorithm}

To measure the goodness of the approximant during the optimization, we introduce the Maximum Absolute Error (MAE) on the validation set defined as follows: 
\begin{equation}
    \label{MAE}
    \mbox{MAE}_{X_{val}, F_{val}}(\textbf{P}_\textbf{f}) = \max_{x_i \in X_{val}, f_i \in F_{val}} | P_f(x_i)- f_i|,
\end{equation} 
where $X_{val}$ and $F_{val}$ are the sets of data locations in the validation set and their corresponding data values and  $\textbf{P}_\textbf{f} = (P_f(x_1), \dots, P_f(x_{k_{val}}))$, with $k_{val} = |X_{val}|$.
In this context, the function to optimize $g$, for which we are searching the maximum, is defined as the negative $\mbox{MAE}_{X_{val}, F_{val}}(\cdot)$ between the known values $F_{val}$ and the approximation evaluated by \textbf{Algorithm} \ref{alg_2} for the points $X_{val}$ and for a specific value of $\varepsilon$. 

\begin{algorithm}
    \caption{$\textbf{RBF}(X, \tilde{X}, F,\Bar{X}, \varepsilon)$}
    \hspace*{\algorithmicindent} \textbf{Input:}
    \vspace{-3mm}
    \begin{itemize}
        \setlength\itemsep{-0.2em}
        \item[] $X$: data locations, $\tilde{X}$: RBF centers, $F$: data values, $\Bar{X}$: evaluation points,  $\varepsilon$: shape parameter.   
    \end{itemize}
    
    \begin{algorithmic}
        \If{$\tilde{X} = X$}
        \State Solve the interpolation system \eqref{linsys}
        \Else 
        \State Solve the approximation system \eqref{linsys_hat}
        \EndIf
    \end{algorithmic}
    \hspace*{\algorithmicindent} \textbf{Output:}
        \vspace{-3mm}
\begin{itemize}
\setlength\itemsep{-0.2em}
    \item[] $\textbf{P}_\textbf{f}$: Evaluation of the approximated or interpolated solution on $\Bar{X}$.
    \end{itemize}
\label{alg_2}
\end{algorithm}

The choice of the next value of $\varepsilon$ to evaluate is driven by an acquisition function $a$ and a surrogate model of $g$ obtained by fitting a Gaussian process. In our case the Expected Improvement (\ref{EI_closed_trade-off}) trades off the exploration and exploitation by means of the $\xi$ parameter (see Subsection \ref{bo}). We underline that the Gaussian process needs some starting iterations to be initialized. Moreover, here we use $niter$ random values of $\varepsilon$ in the search space $\mathcal{X} = (0, \varepsilon_{max}]$ with their associate values of $g$. After the initial settings, the optimization, guided by the acquisition function, performs $niter$ optimization steps, updating the Gaussian process and using it with the acquisition function. Then, the value $\varepsilon^*$ that maximizes the function $g$, namely minimizes the $\mbox{MAE}_{X_{val}, F_{val}}(\cdot)$, is the optimal shape parameter.
Finally, an RBF fit with the found $\varepsilon^*$ is solved and an approximate solution on the points of $\Bar{X}$ is obtained. We refer the reader to \cite{rewiew_BO} for an exhaustive analysis of the BO complexity. A summary of the complete RBF-BO algorithm is sketched in \textbf{Algorithm} \ref{alg_3}.

\begin{algorithm}
\caption{$\textbf{RBF-BO}(X, \tilde{X}, F, \Bar{X}, \mathcal{X}, a, \xi, nstart, niter)$}
\hspace*{\algorithmicindent} \textbf{Input:} 
    \vspace{-3mm}
\begin{itemize}
\setlength\itemsep{-0.2em}
  \item[] $X$: data locations, $\tilde{X}$: RBF centers, 
  $F$: data values,
   $\Bar{X}$: evaluation points, $\mathcal{X}$: parameter search interval, $a$: acquisition function, $\xi$: exploration-exploitation parameter, $nstart$: number of starting iterations,  $niter$: number of Bayesian iterations.
\end{itemize}

\begin{algorithmic}
\State Split $X$, $\tilde{X}$ and $F$ in $X_{train}, X_{val}$,  $\tilde{X}_{train}, \tilde{X}_{val}$, $F_{train}, F_{val}$
\State Set $g \rightarrow -\mbox{MAE}_{X_{val}, F_{val}}(\cdot)$
\State $\varepsilon^* \rightarrow \textbf{BO}(X_{train}, \tilde{X}_{train}, F_{train}, X_{val}, g, \mathcal{X}, a, \xi, nstart, niter)$ \qquad  (call to \textbf{Algorithm} \ref{alg_1})
\State $\textbf{P}_{\textbf{f}} \rightarrow \textbf{RBF}(X, \tilde{X}, F, \Bar{X}, \varepsilon^*)$ \qquad (call to \textbf{Algorithm} \ref{alg_2})
\end{algorithmic} 
\hspace*{\algorithmicindent} \textbf{Output:} 
    \vspace{-3mm}
\begin{itemize}
\setlength\itemsep{-0.2em}
    \item[] $\textbf{P}_\textbf{f}$: Evaluation of the approximated or interpolated solution on $\Bar{X}$.
\end{itemize}
\label{alg_3}
\end{algorithm}

\section{Numerical Experiments}
\label{exp}

In order to apply the Bayesian optimization for radial kernel parameter search, we suppose that the objective function to maximize is 
$-\mbox{MAE}_{X_{val}, F_{val}}(\cdot)$ both in interpolant \eqref{ker_interp} and approximant \eqref{ker_approx} framework. The minus sign ahead is due to the description of BO, in Subsection \ref{bo}, as a maximization process while we are interested in minimization.
All the code was developed in Python 3.9 and the library used to perform the optimization is {\em BayesianOptimization} \cite{BO} in which the default kernel used for the Gaussian process is the  Matérn $5/2$. This section is divided into two subsections: in the former we present some numerical experiments about RBF interpolation, while in the latter we show the algorithm behaviour and performances in a gradual transition from approximation to interpolation. 
To measure the goodness of the approximant, we consider a set of points with the corresponding data values $(\Bar{X}, \Bar{F})$ not used during the fitting of the approximant and  we introduce the MAE on $(\Bar{X}, \Bar{F})$  as follows: 
\begin{equation}
    \mbox{MAE}_{\Bar{X}, \Bar{F}}(\textbf{P}_\textbf{f}) = \max_{x_i \in \Bar{X}, f_i \in \Bar{F}} | P_f(x_i)- f_i|,
\end{equation} 
where   $\textbf{P}_\textbf{f} = (P_f(x_1), \dots, P_f(x_{\bar{k}}))$, with $\bar{k} = |\Bar{X}|$.
All tests are carried out on a MacBook Air (2020), 1.2 GHz Quad-Core Intel Core i7 processor, 16 GB 3733 MHz LPDDR4X RAM,  via Python 3.9.12.

\subsection{Interpolation Results}
\label{experiments_interpolation}
With the aim to make a comparison let us consider the LOOCV to evaluate computation time and error.
The two techniques, LOOCV and BO, are fairly different: on the one hand the LOOCV evaluates the error only on a single point, on the other hand the BO evaluates the error on a set of points.
We use two different data sets in our experiments, the set $(X, F)$ on which we fit the approximant and the set $(\Bar{X}, \Bar{F})$ on which we evaluate the approximation error remaining unchanged during the experiments. 
The use of the set $(X, F)$ depends on the used technique for  $\varepsilon$ detection: recalling that we use Rippa's rule for the evaluation of the LOOCV errors fitting the interpolant on the whole set $X$, we do not perform any subdivision in training and validation sets. 
On the other hand, when BO is used, we perform a further subdivision of $(X, F)$ selecting $80\%$ of the points to have a training set while the remaining part is used as a validation set to compute the error during the optimization. In this way we obtain the subdivision $X_{train}, X_{val}$, $F_{train}, F_{val}$, as previously mentioned in Section \ref{algorithms}.  After finding the corresponding two best values of the shape parameter, we fit a radial kernel interpolant for each on $(X, F)$ and we evaluate the error on $(\bar{X}, \bar{F})$. Acting in this way, the found parameter from LOOCV is favoured because we use the same set of data locations and data values for its determination and for the fitting of the final model. The $\varepsilon$ values among which the LOOCV errors were valuated are $500$ equally spaced values in  $\mathcal{X} = (0,\varepsilon_{max}]=(0, 20]$ and the number of random initial iterations for BO are $5$, while the iteration steps are $25$. Moreover, we control the exploration-exploitation trade-off by means of the $\xi$ parameter (see Subsection  \ref{bo}), considering three values for the parameter:  $\xi = 0.1$ (prefer exploration), $\xi = 0.01$, and $\xi = 0.001$ (prefer exploitation). Moreover, we compare the LOOCV and the BO with the LOOCV$^*$ that consists of a LOOCV with a univariate optimizer. This variant does not try all the $500$ values of $\varepsilon$  but the search of the parameter is driven by the function \em{minimize} of the library \em{scipy.optimize} \cite{scipy} that has been applied in its default form with a starting value set equal to $10$ as the median value of the parameter space search $\mathcal{X}$.
We perform the experiments on random and Halton data locations in the domain $\Omega = [0, 1]^2$ using three different sizes $n$ for each kind of data locations, three different RBFs, i.e.,
\[
\varphi(\varepsilon r) =\left\{
\begin{array}{llll}
\exp(-\varepsilon^2 r^2), & & \quad \mbox{Gaussian $C^{\infty}$} & (\mbox{GA})  \medskip  \\ 
\exp(-\varepsilon r) (\varepsilon r+1), & & \quad \mbox{Mat$\acute{\text{e}}$rn $C^2$} & (\mbox{M2})  \medskip  \\ 
\max \left(1-\varepsilon r,0\right)^4(4\varepsilon r+1), & & \quad \mbox{Wendland $C^2$} & (\mbox{W2}) 
\end{array}
\right.
\]
and the following two test functions \cite{cav19,ren99}
\begin{align*}
\begin{array}{l}
    \displaystyle{f_1(\boldsymbol{x})  =  0.75 \exp\left[{-\frac{(9x_1-2)^2}{4}-\frac{(9x_2-2)^2}{4}}\right] + 0.75 \exp\left[{-\frac{(9x_1-2)^2}{49} - \frac{9x_2+1}{10}}\right] +} \medskip \\  
    \displaystyle{\qquad \quad +  0.5 \exp\left[{-\frac{(9x_1-7)^2}{4}-\frac{(9x_2-3)^2}{4}}\right] -0.2 \exp\left[{-(9x_1-4)^2-(9x_2-7)^2}\right],} \medskip \\
    \displaystyle{f_2(\boldsymbol{x})=  \frac{\sqrt{64 - 81 ((x_1 - 0.5)^2 + (x_2 - 0.5)^2 )}}{9} - 0.5.}
\end{array}
\end{align*}
Notice that the above GA and M2 kernels are globally supported, while the W2 kernel is compactly supported.
The results in Tables \ref{franke_random}--\ref{w2_halton} show, with a few exceptions, that the errors obtained using the LOOCV and BO are pretty similar and the main difference lies in the computational time which is in most cases lower in BO, in particular for high values of $n$ where the cost of inversion of the interpolation matrix is significant. The cases in which the BO does not attain the same precision of LOOCV owe this behaviour to the lower regularity of the function to be approximated and  this can be recovered by increasing the number of Bayesian iterations of the process.  We  also see that a  parameter $\xi = 0.01$ is suitable in most cases. The results obtained by LOOCV$^*$ are for $n = 1000$ and $n = 500$ comparable with the BO results. In some cases the time expenses is less for LOOCV$^*$ but BO has a better error. This could be due to the fact that LOOCV$^*$ gets stuck in a local minima. For $n = 250$ the LOOCV$^*$ behaves better than the BO that is constrained to attain a fixed number of iterations during the exploration of the space. The LOOCV$^*$ method relying on a stop tolerance can save computational expenses. The two methods are comparable but for problems with a high number of data locations BO would be preferred for its speed of convergence and the independence from the starting value. 



\begin{table}[ht!]
\centering
\begin{adjustbox}{width=1\textwidth, center=\textwidth}
\begin{tabular}{|*{9}{p{14mm}|}}
    \cline{4-9}
    \multicolumn{3}{c}{} & \multicolumn{3}{|c}{Matérn kernel} &\multicolumn{3}{|c|}{Gaussian kernel}\\
    \hline
    $n$ &  method & $\xi$ &  time (s) &  $\mbox{MAE}_{\Bar{X}, \Bar{F}}$  & $\varepsilon^*$ & time (s) &  $\mbox{MAE}_{\Bar{X}, \Bar{F}}$  & $\varepsilon^*$ \\
    \hline
    \multirow{5}{*}{1000}  & LOOCV &  & 6.87e+01 &  6.76e-04 &  1.520000 &  7.12e+01 &  8.53e-05 &  7.040000 \\
    \cline{2-9}
    \multirow{5}{*}{}  & LOOCV$^*$ &  & 1.97e+00	& 3.23e-03	& 9.990318 &  5.75e+00	& 4.47e-04	& 10.000001 \\
    \cline{2-9}
    \multirow{5}{*}{} & \multirow{3}{*}{BO} & 0.1 & 5.32e+00 &  6.65e-04 &  1.147200 &  4.98e+00 &  3.40e-04 &  6.159917 \\
    \cline{3-9}
    \multirow{5}{*}{}  & \multirow{3}{*}{} & 0.01  &  6.18e+00 &  6.64e-04 &  1.129313 &  4.41e+00 &  1.66e-03 &  7.042178 \\
    \cline{3-9}
    \multirow{5}{*}{}  & \multirow{3}{*}{} & 0.001 & 4.33e+00 &  6.63e-04 &  0.788350 &  4.33e+00 &  1.45e-04 &  6.465295 \\
    \hline
    \multirow{5}{*}{500}  & LOOCV & &  1.04e+01 &  2.81e-03 &  1.680000 &  1.10e+01 &  7.40e-04 &  5.920000 \\
    \cline{2-9}
    \multirow{5}{*}{}  & LOOCV$^*$ &  & 6.08e-01 &	2.81e-03 & 1.705232 & 4.01e+00 &	1.57e-04	& 6.346624 \\
    \cline{2-9}
    \multirow{5}{*}{}  & \multirow{3}{*}{BO} & 0.1  & 1.90e+00 &  2.73e-02 &  0.001233 &  2.36e+00 &  1.14e-02 &  5.656844 \\
    \cline{3-9}
    \multirow{5}{*}{}  & \multirow{3}{*}{} & 0.01 &  1.72e+00 &  2.83e-03 &  1.453144 &  2.45e+00 &  1.99e-01 &  5.663782 \\
    \cline{3-9}
    \multirow{5}{*}{}  & \multirow{3}{*}{} & 0.001 & 1.98e+00 &  2.60e-02 &  0.001000 &  3.30e+00 &  1.10e-03 &  5.906976 \\
    \hline
    \multirow{5}{*}{250}  & LOOCV & &2.23e+00 &  4.50e-03 &  0.840000 &  2.11e+00 &  4.74e-02 &  5.520000 \\
    \cline{2-9}
    \multirow{5}{*}{}  & LOOCV$^*$ &  & 7.60e-01&	4.51e-03&	0.818112 &  5.33e-01 &	3.72e-02 &	6.293530 \\
    \cline{2-9}
    \multirow{5}{*}{}  & \multirow{3}{*}{BO} & 0.1  &  1.52e+00 &  4.35e-03 &  1.147200 &  1.80e+00 &  5.00e-02 &  5.475028 \\
    \cline{3-9}
    \multirow{5}{*}{}  & \multirow{3}{*}{} & 0.01 & 1.48e+00 &  4.20e-03 &  1.528790 &  1.86e+00 &  3.51e-02 &  6.069408 \\
    \cline{3-9}
    \multirow{5}{*}{}  & \multirow{3}{*}{} & 0.001 & 1.54e+00 &  4.19e-03 &  1.570866 &  2.55e+00 &  4.44e-02 &  6.630227 \\
    \hline
\end{tabular}
\end{adjustbox}
\caption{Computational time, $\mbox{MAE}_{\Bar{X}, \Bar{F}}$ and shape parameter using M2 and GA kernels for the interpolation of $f_1$ on various sets of random points by LOOCV, LOOCV$^*$ and BO.}
\label{franke_random}
\end{table}

\begin{table}[ht!]
\centering
\begin{adjustbox}{width=1\textwidth, center=\textwidth}
\begin{tabular}{|*{9}{p{14mm}|}}
    \cline{4-9}
    \multicolumn{3}{c}{} & \multicolumn{3}{|c}{Matérn kernel} &\multicolumn{3}{|c|}{Gaussian kernel}\\
    \hline
         $n$ &  method & $\xi$ &  time (s) &  $\mbox{MAE}_{\Bar{X}, \Bar{F}}$  & $\varepsilon^*$ & time (s) &  $\mbox{MAE}_{\Bar{X}, \Bar{F}}$  & $\varepsilon^*$ \\
    \hline
    \multirow{5}{*}{1000}  & LOOCV &  & 6.09e+01 &  7.69e-04 &  1.240000 &  6.60e+01 &  1.29e-05 &  6.440000 \\
    \cline{2-9}
    \multirow{5}{*}{}  & LOOCV$^*$ &  & 3.87e+00&	4.05e-04&	1.251510 &  6.11e+00	&3.72e-04	&10.000009 \\
    \cline{2-9}
    \multirow{5}{*}{} & \multirow{3}{*}{BO} & 0.1 & 4.72e+00 &  7.70e-04 &  1.147200 &  4.42e+00 &  2.91e-05 &  6.153352 \\
    \cline{3-9}
    \multirow{5}{*}{}  & \multirow{3}{*}{} & 0.01  & 6.06e+00 &  7.68e-04 &  1.604766 &  4.31e+00 &  2.93e-05 &  6.847236 \\
    \cline{3-9}
    \multirow{5}{*}{}  & \multirow{3}{*}{} & 0.001 & 6.20e+00 &  7.68e-04 &  1.328962 &  4.29e+00 &  9.71e-06 &  7.491428 \\
    \hline
    \multirow{5}{*}{500}  & LOOCV & &   7.59e+00 &  3.33e-03 &  4.640000 &  7.96e+00 &  1.08e-04 &  6.080000 \\
    \cline{2-9}
    \multirow{5}{*}{}  & LOOCV$^*$ &  & 5.13e-01 &	2.63e-03 &	1.495576 &  1.81e+00&	1.56e-03&	6.782924 \\
    \cline{2-9}
    \multirow{5}{*}{}  & \multirow{3}{*}{BO} & 0.1  &  1.54e+00 &  2.32e-03 &  1.147200 &  2.81e+00 &  1.06e-04 &  6.064566 \\
    \cline{3-9}
    \multirow{5}{*}{}  & \multirow{3}{*}{} & 0.01 &  1.58e+00 &  2.32e-03 &  1.147200 &  2.55e+00 &  1.61e-04 &  6.116520 \\
    \cline{3-9}
    \multirow{5}{*}{}  & \multirow{3}{*}{} & 0.001 &  1.61e+00 &  3.32e-03 &  0.001000 &  2.19e+00 &  1.20e-04 &  6.075517 \\
    \hline
    \multirow{5}{*}{250}  & LOOCV & & 1.24e+00 &  6.52e-03 &  2.280000 &  1.21e+00 &  5.82e-03 &  5.800000 \\
    \cline{2-9}
    \multirow{5}{*}{}  & LOOCV$^*$ &  & 3.82e-01	& 9.87e-03	& 3.365488 &  5.26e-01	& 2.58e-02	& 4.997309 \\
    \cline{2-9}
    \multirow{5}{*}{}  & \multirow{3}{*}{BO} & 0.1  &   1.49e+00 &  7.61e-03 &  3.639104 &  1.63e+00 &  6.31e-03 &  5.957693 \\
    \cline{3-9}
    \multirow{5}{*}{}  & \multirow{3}{*}{} & 0.01 &  1.49e+00 &  7.61e-03 &  3.639104 &  2.12e+00 &  5.93e-03 &  5.844204 \\
    \cline{3-9}
    \multirow{5}{*}{}  & \multirow{3}{*}{} & 0.001 & 1.39e+00 &  8.06e-03 &  3.981307 &  2.59e+00 &  5.90e-03 &  5.834172 \\
    \hline
\end{tabular}
\end{adjustbox}
\caption{Computational time, $\mbox{MAE}_{\Bar{X}, \Bar{F}}$ and shape parameter using M2 and GA kernels for the interpolation of $f_1$ on various sets of Halton points by LOOCV, LOOCV$^*$ and BO.}
\label{franke_halton}
\end{table}

\begin{table}[ht!]
\centering
\begin{adjustbox}{width=1\textwidth, center=\textwidth}
\begin{tabular}{|*{9}{p{14mm}|}}
    \cline{4-9}
    \multicolumn{3}{c}{} & \multicolumn{3}{|c}{Matérn kernel} &\multicolumn{3}{|c|}{Gaussian kernel}\\
    \hline
         $n$ &  method & $\xi$ &  time (s) &  $\mbox{MAE}_{\Bar{X}, \Bar{F}}$  & $\varepsilon^*$ & time (s) &  $\mbox{MAE}_{\Bar{X}, \Bar{F}}$  & $\varepsilon^*$ \\
    \hline
    \multirow{5}{*}{1000}  & LOOCV &  & 7.05e+01 &  8.60e-04 &  0.040000 &  7.90e+01 &  2.66e-05 &  6.000000 \\
    \cline{2-9}
    \multirow{5}{*}{}  & LOOCV$^*$ &  & 3.73e+00 &	1.18e-03&	0.668190 &  7.16e+00 &	2.91e-04	& 10.000001 \\
    \cline{2-9}
    \multirow{5}{*}{} & \multirow{3}{*}{BO} & 0.1 & 3.67e+00 &  1.41e-03 &  1.147200 &  4.28e+00 &  1.50e-05 &  4.349624 \\
    \cline{3-9}
    \multirow{5}{*}{}  & \multirow{3}{*}{} & 0.01  & 3.69e+00 &  1.59e-03 &  1.507068 &  4.62e+00 &  7.28e-04 &  3.398101 \\
    \cline{3-9}
    \multirow{5}{*}{}  & \multirow{3}{*}{} & 0.001 & 3.70e+00 &  1.08e-03 &  0.474564 &  5.20e+00 &  2.07e-05 &  4.782848 \\
    \hline
    \multirow{5}{*}{500}  & LOOCV & &    1.11e+01 &  2.24e-03 &  0.040000 &  1.10e+01 &  4.14e-05 &  3.840000 \\
    \cline{2-9}
    \multirow{5}{*}{}  & LOOCV$^*$ &  & 8.44e-01 &	2.75e-03 &	0.421860 &  1.96e+00 &	1.15e-03 &	7.824069 \\
    \cline{2-9}
    \multirow{5}{*}{}  & \multirow{3}{*}{BO} & 0.1  &  1.78e+00 &  8.53e-03 &  0.001233 &  1.77e+00 &  7.23e-06 &  3.121217 \\
    \cline{3-9}
    \multirow{5}{*}{}  & \multirow{3}{*}{} & 0.01 &   1.77e+00 &  8.53e-03 &  0.001233 &  2.13e+00 &  7.23e-06 &  3.121217 \\
    \cline{3-9}
    \multirow{5}{*}{}  & \multirow{3}{*}{} & 0.001 &   1.65e+00 &  2.71e-03 &  0.393487 &  2.66e+00 &  7.23e-06 &  3.121217 \\
    \hline
    \multirow{5}{*}{250}  & LOOCV & & 2.09e+00 &  1.20e-02 &  0.040000 &  2.09e+00 &  2.17e-04 &  2.760000 \\
    \cline{2-9}
    \multirow{5}{*}{}  & LOOCV$^*$ &  & 2.48e-01	& 1.27e-02&	0.143446 &  4.90e-01 &	7.55e-04 &	1.136128 \\
    \cline{2-9}
    \multirow{5}{*}{}  & \multirow{3}{*}{BO} & 0.1  &  1.32e+00 &  3.07e-02 &  7.174972 &  2.53e+00 &  3.81e-04 &  1.867692 \\
    \cline{3-9}
    \multirow{5}{*}{}  & \multirow{3}{*}{} & 0.01 & 1.40e+00 &  3.07e-02 &  7.174972 &  1.77e+00 &  9.57e-05 &  2.728965 \\
    \cline{3-9}
    \multirow{5}{*}{}  & \multirow{3}{*}{} & 0.001 &  1.44e+00 &  1.13e-02 &  0.001000 &  2.04e+00 &  4.43e-04 &  2.979949 \\
    \hline
\end{tabular}
\end{adjustbox}
\caption{Computational time, $\mbox{MAE}_{\Bar{X}, \Bar{F}}$ and shape parameter using M2 and GA kernels for the interpolation of $f_2$ on various sets of random points by LOOCV, LOOCV$^*$ and BO.}
\label{f6_random}
\end{table}

\begin{table}[ht!]
\centering
\begin{adjustbox}{width=1\textwidth, center=\textwidth}
\begin{tabular}{|*{9}{p{14mm}|}}
    \cline{4-9}
    \multicolumn{3}{c}{} & \multicolumn{3}{|c}{Matérn kernel} &\multicolumn{3}{|c|}{Gaussian kernel}\\
    \hline
        $n$ &  method & $\xi$ &  time (s) &  $\mbox{MAE}_{\Bar{X}, \Bar{F}}$  & $\varepsilon^*$ & time (s) &  $\mbox{MAE}_{\Bar{X}, \Bar{F}}$  & $\varepsilon^*$ \\
    \hline
    \multirow{5}{*}{1000}  & LOOCV &  & 5.23e+01 &  2.61e-04 &  0.040000 &  5.07e+01 &  1.44e-06 &  4.320000 \\
    \cline{2-9}
    \multirow{5}{*}{}  & LOOCV$^*$ &  & 3.84e+00	&3.35e-04&	0.257555&  6.11e+00	& 1.21e-04	& 9.999999 \\
    \cline{2-9}
    \multirow{5}{*}{} & \multirow{3}{*}{BO} & 0.1 & 3.10e+00 &  4.48e-04 &  1.147200 &  3.19e+00 &  1.79e-06 &  3.121217 \\
    \cline{3-9}
    \multirow{5}{*}{}  & \multirow{3}{*}{} & 0.01  & 3.30e+00 &  3.94e-04 &  0.838285 &  3.68e+00 &  1.91e-06 &  3.402787 \\
    \cline{3-9}
    \multirow{5}{*}{}  & \multirow{3}{*}{} & 0.001 & 3.37e+00 &  4.11e-04 &  0.934542 &  3.87e+00 &  1.79e-06 &  3.121217 \\
    \hline
    \multirow{5}{*}{500}  & LOOCV & &    8.49e+00 &  7.10e-04 &  0.040000 &  7.46e+00 &  2.10e-05 &  4.840000 \\
    \cline{2-9}
    \multirow{5}{*}{}  & LOOCV$^*$ &  & 6.22e-02	& 7.01e-03&	9.999965 &  1.17e+00 &	3.06e-04	& 5.917107 \\
    \cline{2-9}
    \multirow{5}{*}{}  & \multirow{3}{*}{BO} & 0.1  & 1.51e+00 &  1.12e-03 &  1.147200 &  1.69e+00 &  7.76e-06 &  2.199642 \\
    \cline{3-9}
    \multirow{5}{*}{}  & \multirow{3}{*}{} & 0.01 &   1.51e+00 &  6.61e-04 &  0.001233 &  2.54e+00 &  1.70e-05 &  3.378652 \\
    \cline{3-9}
    \multirow{5}{*}{}  & \multirow{3}{*}{} & 0.001 &  1.52e+00 &  1.22e-03 &  0.001000 &  2.66e+00 &  1.43e-05 &  3.377515 \\
    \hline
    \multirow{5}{*}{250}  & LOOCV & & 1.23e+00 &  3.73e-03 &  0.040000 &  1.20e+00 &  4.09e-05 &  2.240000 \\
    \cline{2-9}
    \multirow{5}{*}{}  & LOOCV$^*$ &  & 4.21e-01	& 2.43e-03	& 0.020730 &  5.96e-01	& 6.39e-04	& 4.251622 \\
    \cline{2-9}
    \multirow{5}{*}{}  & \multirow{3}{*}{BO} & 0.1  &   1.37e+00 &  3.64e-03 &  0.001233 &  2.74e+00 &  1.70e-03 &  1.414414 \\
    \cline{3-9}
    \multirow{5}{*}{}  & \multirow{3}{*}{} & 0.01 & 1.39e+00 &  3.64e-03 &  0.001233 &  2.32e+00 &  3.03e-05 &  2.057707 \\
    \cline{3-9}
    \multirow{5}{*}{}  & \multirow{3}{*}{} & 0.001 &  1.34e+00 &  8.87e-03 &  3.121217 &  1.67e+00 &  4.20e-05 &  2.820584 \\
    \hline
\end{tabular}
\end{adjustbox}
\caption{Computational time, $\mbox{MAE}_{\Bar{X}, \Bar{F}}$ and shape parameter using M2 and GA kernels for the interpolation of $f_2$ on various sets of Halton points by LOOCV, LOOCV$^*$ and BO.}
\label{f6_halton}
\end{table}

\begin{table}[ht!]
\centering
\begin{adjustbox}{width=1\textwidth, center=\textwidth}
\begin{tabular}{|*{9}{p{14mm}|}}
    \cline{4-9}
    \multicolumn{3}{c}{} & \multicolumn{3}{|c}{$f_1$} &\multicolumn{3}{|c|}{$f_2$}\\
    \hline
       $n$ &  method & $\xi$ &  time (s) &  $\mbox{MAE}_{\Bar{X}, \Bar{F}}$  & $\varepsilon^*$ & time (s) &  $\mbox{MAE}_{\Bar{X}, \Bar{F}}$  & $\varepsilon^*$ \\
    \hline
    \multirow{5}{*}{1000}  & LOOCV &  &  9.37e+01 &  6.77e-04 &  0.280000 &  7.77e+01 &  9.13e-04 &  0.040000 \\
    \cline{2-9}
    \multirow{5}{*}{}  & LOOCV$^*$ &  & 1.33e+01&	6.71e-04&	0.240126 &  8.86e+00	&8.95e-04	&0.030410 \\
    \cline{2-9}
    \multirow{5}{*}{} & \multirow{3}{*}{BO} & 0.1 & 4.86e+00 &  7.37e-04 &  0.001000 &  5.35e+00 &  8.34e-04 &  0.001000 \\
    \cline{3-9}
    \multirow{5}{*}{}  & \multirow{3}{*}{} & 0.01  & 5.87e+00 &  6.85e-04 &  0.001000 &  4.97e+00 &  8.30e-04 &  0.001000 \\
    \cline{3-9}
    \multirow{5}{*}{}  & \multirow{3}{*}{} & 0.001 &  6.07e+00 &  6.85e-04 &  0.001000 &  6.59e+00 &  8.30e-04 &  0.001000 \\
    \hline
    \multirow{5}{*}{500}  & LOOCV & &    1.25e+01 &  2.95e-03 &  0.280000 &  1.24e+01 &  2.38e-03 &  0.040000 \\
    \cline{2-9}
    \multirow{5}{*}{}  & LOOCV$^*$ &  & 1.59e+00	& 2.95e-03&	0.280290 &  1.66e+00 &	5.08e-03 &	0.818422 \\
    \cline{2-9}
    \multirow{5}{*}{}  & \multirow{3}{*}{BO} & 0.1  &  1.86e+00 &  3.44e-03 &  0.001000 &  1.82e+00 &  2.18e-03 &  0.001000 \\
    \cline{3-9}
    \multirow{5}{*}{}  & \multirow{3}{*}{} & 0.01 &  2.09e+00 &  3.49e-03 &  0.001000 &  1.87e+00 &  2.19e-03 &  0.001000 \\
    \cline{3-9}
    \multirow{5}{*}{}  & \multirow{3}{*}{} & 0.001 & 2.49e+00 &  3.49e-03 &  0.001000 &  2.33e+00 &  2.21e-03 &  0.001096 \\
    \hline
    \multirow{5}{*}{250}  & LOOCV & &  2.82e+00 &  4.69e-03 &  0.320000 &  2.53e+00 &  1.27e-02 &  0.040000 \\
    \cline{2-9}
    \multirow{5}{*}{}  & LOOCV$^*$ &  & 1.63e+00	& 4.68e-03	& 0.287048 &  3.15e-01 &	1.27e-02 &	0.037173 \\
    \cline{2-9}
    \multirow{5}{*}{}  & \multirow{3}{*}{BO} & 0.1  &  1.56e+00 &  4.82e-03 &  0.440938 &  1.54e+00 &  2.34e-02 &  1.249925 \\
    \cline{3-9}
    \multirow{5}{*}{}  & \multirow{3}{*}{} & 0.01 &   1.56e+00 &  4.68e-03 &  0.280038 &  1.41e+00 &  2.34e-02 &  1.237074 \\
    \cline{3-9}
    \multirow{5}{*}{}  & \multirow{3}{*}{} & 0.001 &  2.41e+00 &  4.68e-03 &  0.270527 &  1.61e+00 &  2.33e-02 &  1.262232 \\
    \hline
\end{tabular}
\end{adjustbox}
\caption{Computational time, $\mbox{MAE}_{\Bar{X}, \Bar{F}}$ and shape parameter using W2  kernel for the interpolation of $f_1$  and $f_2$ on various sets of random points by LOOCV, LOOCV$^*$ and BO.}
\label{w2_random}
\end{table}

\begin{table}[ht]
\centering
\begin{adjustbox}{width=1\textwidth, center=\textwidth}
\begin{tabular}{|*{9}{p{14mm}|}}
    \cline{4-9}
    \multicolumn{3}{c}{} & \multicolumn{3}{|c}{$f_1$} &\multicolumn{3}{|c|}{$f_2$}\\
    \hline
       $n$ &  method & $\xi$ &  time (s) &  $\mbox{MAE}_{\Bar{X}, \Bar{F}}$  & $\varepsilon^*$ & time (s) &  $\mbox{MAE}_{\Bar{X}, \Bar{F}}$  & $\varepsilon^*$ \\
    \hline
    \multirow{5}{*}{1000}  & LOOCV &  &   5.86e+01 &  7.91e-04 &  0.240000 &  5.96e+01 &  2.76e-04 &  0.040000 \\
    \cline{2-9}
    \multirow{5}{*}{}  & LOOCV$^*$ &  & 6.14e+00 &	4.58e-04 &	0.179607 &  1.84e+01	& 3.21e-04 &	0.046904 \\
    \cline{2-9}
    \multirow{5}{*}{} & \multirow{3}{*}{BO} & 0.1 &3.59e+00 &  8.25e-04 &  0.001000 &  3.25e+00 &  2.56e-04 &  0.001233 \\
    \cline{3-9}
    \multirow{5}{*}{}  & \multirow{3}{*}{} & 0.01  & 3.61e+00 &  7.91e-04 &  0.243379 &  3.29e+00 &  2.57e-04 &  0.001000 \\
    \cline{3-9}
    \multirow{5}{*}{}  & \multirow{3}{*}{} & 0.001 &   4.12e+00 &  7.93e-04 &  0.277008 &  4.20e+00 &  2.55e-04 &  0.001096 \\
    \hline
    \multirow{5}{*}{500}  & LOOCV & &  8.77e+00 &  3.81e-03 &  1.000000 &  8.88e+00 &  7.48e-04 &  0.040000 \\
    \cline{2-9}
    \multirow{5}{*}{}  & LOOCV$^*$ &  & 2.07e+00 &	3.09e-03 &	0.263759&  2.29e+00	& 1.69e-03&	0.019589 \\
    \cline{2-9}
    \multirow{5}{*}{}  & \multirow{3}{*}{BO} & 0.1  &  1.65e+00 &  2.53e-03 &  0.001000 &  1.69e+00 &  6.96e-04 &  0.001000 \\
    \cline{3-9}
    \multirow{5}{*}{}  & \multirow{3}{*}{} & 0.01 & 1.84e+00 &  2.53e-03 &  0.001792 &  1.71e+00 &  6.96e-04 &  0.001000 \\
    \cline{3-9}
    \multirow{5}{*}{}  & \multirow{3}{*}{} & 0.001 & 2.18e+00 &  2.53e-03 &  0.001096 &  2.13e+00 &  6.97e-04 &  0.001096 \\
    \hline
    \multirow{5}{*}{250}  & LOOCV & &  1.58e+00 &  6.67e-03 &  0.360000 &  1.53e+00 &  3.94e-03 &  0.040000 \\
    \cline{2-9}
    \multirow{5}{*}{}  & LOOCV$^*$ &  & 2.46e-01	& 1.18e-02	& 1.064206 &  9.70e-01	& 2.46e-03 &	0.011255 \\
    \cline{2-9}
    \multirow{5}{*}{}  & \multirow{3}{*}{BO} & 0.1  & 1.41e+00 &  7.18e-03 &  0.568376 &  1.47e+00 &  3.65e-03 &  0.001000 \\
    \cline{3-9}
    \multirow{5}{*}{}  & \multirow{3}{*}{} & 0.01 &   1.99e+00 &  6.81e-03 &  0.449130 &  1.44e+00 &  3.65e-03 &  0.001000 \\
    \cline{3-9}
    \multirow{5}{*}{}  & \multirow{3}{*}{} & 0.001 & 1.94e+00 &  6.70e-03 &  0.379231 &  1.84e+00 &  3.65e-03 &  0.001000 \\
    \hline
\end{tabular}
\end{adjustbox}
\caption{Computational time, $\mbox{MAE}_{\Bar{X}, \Bar{F}}$ and shape parameter using W2 kernel for the interpolation of $f_1$  and $f_2$ on various sets of Halton points by LOOCV, LOOCV$^*$ and BO.}
\label{w2_halton}
\end{table}

\subsection{Approximation Results}
In this subsection we show some experiments of approximation by varying the percentage of data locations in $X$ used as approximation centers. We note that in the case of $100\%$ of nodes used as centers we fall into the interpolation case and thus the results are pretty similar (not exactly the same due to the stochasticity of the experiments) to the ones obtained in Subsection \ref{experiments_interpolation}. In all cases, considering approximation instead of interpolation leads to a loss of accuracy with the benefit of time savings. Notably, in  Tables \ref{tab:approx_gaussian}--\ref{tab:approx_wendland} we can see three different situations: 
\begin{enumerate}
\item When the training set is small, as in the case of $n=250$ points, the precision is low in all cases and there is nearly no saving time in considering approximation, thus, in this case, interpolation is recommended to preserve accuracy.
\item When the training set is medium-size, for instance $n=500$ points, we can see that in the majority of cases an approximation with $80\%$ of points as centers saves computational time with  respect to the interpolation and it is therefore to be preferred. 
\item When the training set is large-size, i.e. $n=1000$ points, we can see that in all cases we have a more significant decrease in computational time if we perform approximation using  $80\%$ of points as centers with respect to the interpolation, without compromising precision, that can motivate in a strong way the choice of approximation.
\end{enumerate}
In conclusion, with larger datasets the saving in computational time assumes a greater magnitude and could be a fair trade-off between computational expense and accuracy. 
\begin{table}[ht]
\centering
\begin{adjustbox}{width=1\textwidth, center=\textwidth}
\begin{tabular}{|*{8}{p{14mm}|}}
    \cline{3-8}
    \multicolumn{2}{c}{} & \multicolumn{3}{|c}{Random} &\multicolumn{3}{|c|}{Halton}\\
    \hline
       $n$  & centers($ \%$) &  time (s) &  $\mbox{MAE}_{\Bar{X}, \Bar{F}}$  & $\varepsilon^*$ & time (s) &  $\mbox{MAE}_{\Bar{X}, \Bar{F}}$  & $\varepsilon^*$ \\
    \hline
    \multirow{2}{*}{1000} &   20 &  3.79e+00 &  2.96e-03 &  4.861528 &  2.98e+00 &  1.88e-03 &  5.147405 \\
    \cline{2-8}
    \multirow{2}{*}{} & 40 &  5.59e+00 &  1.50e-03 &  5.285210 &  2.90e+00 &  1.05e-04 &  5.594780 \\
    \cline{2-8}
    \multirow{2}{*}{}  & 60 &  5.46e+00 &  2.33e-04 &  5.744282 &  3.92e+00 &  2.78e-05 &  6.348615 \\
    \cline{2-8}
    \multirow{2}{*}{}  & 80 &  5.85e+00 &  3.90e-05 &  6.319368 &  5.45e+00 &  2.43e-05 &  5.698942 \\
    \cline{2-8}
    \multirow{2}{*}{}  & 100 &  7.61e+00 &  7.87e-05 &  6.319614 &  6.35e+00 &  9.87e-06 &  6.309349 \\
    \hline
    \multirow{2}{*}{500}  & 20 &  1.64e+00 &  1.98e-02 &  4.471739 &  1.57e+00 &  1.67e-02 &  5.361774 \\
    \cline{2-8}
    \multirow{2}{*}{}  & 40 &  2.48e+00 &  6.83e-03 &  5.289607 &  1.83e+00 &  1.68e-03 &  5.171920 \\
    \cline{2-8}
    \multirow{2}{*}{}  & 60 &  2.82e+00 &  2.89e-03 &  5.869333 &  2.33e+00 &  5.12e-04 &  5.336932 \\
    \cline{2-8}
    \multirow{2}{*}{}  & 80 &  2.60e+00 &  5.61e-03 &  5.043093 &  2.07e+00 &  6.12e-04 &  5.312729 \\
    \cline{2-8}
    \multirow{2}{*}{}  & 100 &  4.33e+00 &  2.34e-03 &  6.297596 &  2.20e+00 &  5.31e-04 &  5.459864 \\
    \hline
    \multirow{2}{*}{250}  & 20 &  1.35e+00 &  7.73e-02 &  3.121217 &  1.43e+00 &  5.63e-02 &  3.772420 \\
    \cline{2-8}
    \multirow{2}{*}{}  & 40 &  1.38e+00 &  2.08e-02 &  4.491812 &  1.32e+00 &  1.87e-01 &  1.694190 \\
    \cline{2-8}
    \multirow{2}{*}{}  & 60 &  1.38e+00 &  2.70e-02 &  5.093333 &  1.28e+00 &  1.28e-02 &  4.491812 \\
    \cline{2-8}
    \multirow{2}{*}{}  & 80 &  1.55e+00 &  6.55e-02 &  4.548724 &  2.23e+00 &  1.04e-02 &  5.183731 \\
    \cline{2-8}
    \multirow{2}{*}{}  & 100 &  1.45e+00 &  1.34e-02 &  5.582551 &  1.36e+00 &  6.86e-03 &  5.143759 \\
    \hline
\end{tabular}
\end{adjustbox}
\caption{Computational time, $\mbox{MAE}_{\Bar{X}, \Bar{F}}$ and shape parameter using GA kernel for the approximation of $f_1$ on various sets of random and Halton points as the center percentage  varies by applying LOOCV and BO.}
\label{tab:approx_gaussian}
\end{table}

\begin{table}[ht]
\centering
\begin{adjustbox}{width=1\textwidth, center=\textwidth}
\begin{tabular}{|*{8}{p{14mm}|}}
    \cline{3-8}
    \multicolumn{2}{c}{} & \multicolumn{3}{|c}{Random} &\multicolumn{3}{|c|}{Halton}\\
    \hline
       $n$  & centers($\%$)  &  time (s) &  $\mbox{MAE}_{\Bar{X}, \Bar{F}}$  & $\varepsilon^*$ & time (s) &  $\mbox{MAE}_{\Bar{X}, \Bar{F}}$  & $\varepsilon^*$ \\
    \hline
    \multirow{2}{*}{1000} &   20 &  3.79e+00 &  2.96e-03 &  4.861528 &  2.98e+00 &  1.88e-03 &  5.147405 \\
    \cline{2-8}
    \multirow{2}{*}{} & 40 & 3.21e+00 &  2.01e-03 &  1.498876 &  2.95e+00 &  2.75e-03 &  2.164716 \\
    \cline{2-8}
    \multirow{2}{*}{}  & 60 &   3.60e+00 &  1.74e-03 &  1.312339 &  3.79e+00 &  1.51e-03 &  2.244381 \\
    \cline{2-8}
    \multirow{2}{*}{}  & 80 &   5.22e+00 &  1.02e-03 &  1.332212 &  5.78e+00 &  9.02e-04 &  1.676676 \\
    \cline{2-8}
    \multirow{2}{*}{}  & 100 & 6.41e+00 &  1.03e-03 &  1.332212 &  7.83e+00 &  8.61e-04 &  1.695884 \\
    \hline
    \multirow{2}{*}{500}  & 20 &   3.49e+00 &  3.29e-02 &  7.491428 &  2.61e+00 &  1.64e-02 &  2.245995 \\
    \cline{2-8}
    \multirow{2}{*}{}  & 40 &  2.35e+00 &  1.16e-02 &  1.251539 &  2.96e+00 &  4.82e-03 &  2.087647 \\
    \cline{2-8}
    \multirow{2}{*}{}  & 60 &  3.77e+00 &  1.08e-02 &  2.086028 &  4.79e+00 &  5.07e-03 &  1.086503 \\
    \cline{2-8}
    \multirow{2}{*}{}  & 80 &  2.43e+00 &  4.36e-03 &  1.293058 &  4.43e+00 &  2.88e-03 &  2.647202 \\
    \cline{2-8}
    \multirow{2}{*}{}  & 100 &  2.85e+00 &  4.19e-03 &  1.398220 &  5.37e+00 &  2.51e-03 &  2.183219 \\
    \hline
    \multirow{2}{*}{250}  & 20 &  1.34e+00 &  9.87e-02 &  0.488346 &  1.45e+00 &  4.99e-02 &  5.536034 \\
    \cline{2-8}
    \multirow{2}{*}{}  & 40 & 1.52e+00 &  6.54e-02 &  6.125847 &  1.26e+00 &  1.77e-02 &  2.069197 \\
    \cline{2-8}
    \multirow{2}{*}{}  & 60 &   3.18e+00 &  2.39e-02 &  2.085691 &  1.32e+00 &  9.46e-03 &  2.143143 \\
    \cline{2-8}
    \multirow{2}{*}{}  & 80 &   3.77e+00 &  1.66e-02 &  0.898543 &  2.51e+00 &  4.68e-03 &  1.521938 \\
    \cline{2-8}
    \multirow{2}{*}{}  & 100 &  2.86e+00 &  1.34e-02 &  4.262523 &  2.41e+00 &  4.44e-03 &  1.556272 \\
    \hline
\end{tabular}
\end{adjustbox}
\caption{Computational time, $\mbox{MAE}_{\Bar{X}, \Bar{F}}$ and shape parameter using M2 kernel for the approximation of $f_1$ on various sets of random and Halton points as the center percentage  varies by applying LOOCV and BO.}
\label{tab:approx_matern}
\end{table}

\begin{table}[ht]
\centering
\begin{adjustbox}{width=1\textwidth, center=\textwidth}
\begin{tabular}{|*{8}{p{14mm}|}}
    \cline{3-8}
    \multicolumn{2}{c}{} & \multicolumn{3}{|c}{Random} &\multicolumn{3}{|c|}{Halton}\\
    \hline
       $n$  & centers($\%$) &  time (s) &  $\mbox{MAE}_{\Bar{X}, \Bar{F}}$  & $\varepsilon^*$ & time (s) &  $\mbox{MAE}_{\Bar{X}, \Bar{F}}$  & $\varepsilon^*$ \\
    \hline
    \multirow{2}{*}{1000} &   20 &  3.89e+00 &  8.69e-03 &  0.293414 &  1.80e+00 &  1.30e-02 &  0.525864 \\
    \cline{2-8}
    \multirow{2}{*}{} & 40 &   3.84e+00 &  2.16e-03 &  0.279970 &  3.35e+00 &  6.61e-03 &  0.382713 \\
    \cline{2-8}
    \multirow{2}{*}{}  & 60 &  5.80e+00 &  6.83e-03 &  0.286352 &  3.89e+00 &  3.89e-03 &  0.294302 \\
    \cline{2-8}
    \multirow{2}{*}{}  & 80 &  7.76e+00 &  6.66e-03 &  0.266835 &  5.52e+00 &  8.00e-04 &  0.320220 \\
    \cline{2-8}
    \multirow{2}{*}{}  & 100 & 8.81e+00 &  1.02e-03 &  0.310954 &  7.05e+00 &  8.59e-04 &  0.328610 \\
    \hline
    \multirow{2}{*}{500}  & 20 &  2.23e+00 &  2.49e-02 &  0.127089 &  2.18e+00 &  3.22e-02 &  0.344410 \\
    \cline{2-8}
    \multirow{2}{*}{}  & 40 & 2.59e+00 &  8.41e-03 &  0.199470 &  2.14e+00 &  5.58e-03 &  0.337803 \\
    \cline{2-8}
    \multirow{2}{*}{}  & 60 &  3.40e+00 &  8.91e-03 &  0.323238 &  2.15e+00 &  4.47e-03 &  0.299074 \\
    \cline{2-8}
    \multirow{2}{*}{}  & 80 &  2.64e+00 &  6.81e-03 &  1.135019 &  2.42e+00 &  2.68e-03 &  0.261940 \\
    \cline{2-8}
    \multirow{2}{*}{}  & 100 &  2.98e+00 &  6.74e-03 &  1.120099 &  2.10e+00 &  2.60e-03 &  0.250589 \\
    \hline
    \multirow{2}{*}{250}  & 20 &  1.26e+00 &  7.70e-02 &  1.065228 &  1.38e+00 &  7.97e-02 &  1.714649 \\
    \cline{2-8}
    \multirow{2}{*}{}  & 40 & 1.51e+00 &  4.42e-02 &  0.311686 &  1.54e+00 &  3.12e-02 &  0.368420 \\
    \cline{2-8}
    \multirow{2}{*}{}  & 60 & 1.53e+00 &  2.38e-02 &  0.736349 &  1.41e+00 &  1.67e-02 &  0.401290 \\
    \cline{2-8}
    \multirow{2}{*}{}  & 80 & 1.61e+00 &  1.76e-02 &  0.425162 &  1.50e+00 &  7.12e-03 &  0.311549 \\
    \cline{2-8}
    \multirow{2}{*}{}  & 100 &  1.74e+00 &  1.60e-02 &  0.423409 &  1.72e+00 &  4.63e-03 &  0.281553 \\
    \hline
\end{tabular}
\end{adjustbox}
\caption{Computational time, $\mbox{MAE}_{\Bar{X}, \Bar{F}}$ and shape parameter using W2 kernel for the approximation of $f_1$ on various sets of random and Halton points as the center percentage  varies by applying LOOCV and BO.}
\label{tab:approx_wendland}
\end{table}

\section{Application to Real Data}
\label{real_data_example}
\begin{figure}[ht]
    \centering
    \includegraphics[scale = 0.5]{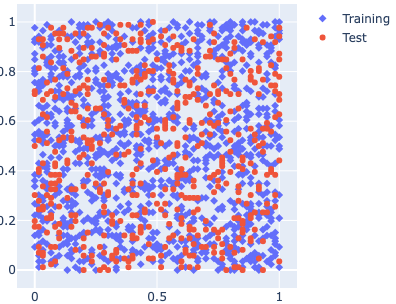}
    \caption{Training (blue) and test (red) sets extracted from the \emph{volcano} dataset.}
    \label{fig:train_test}
\end{figure}

\begin{figure}[ht]
    \centering
    \includegraphics[width=.49\textwidth]{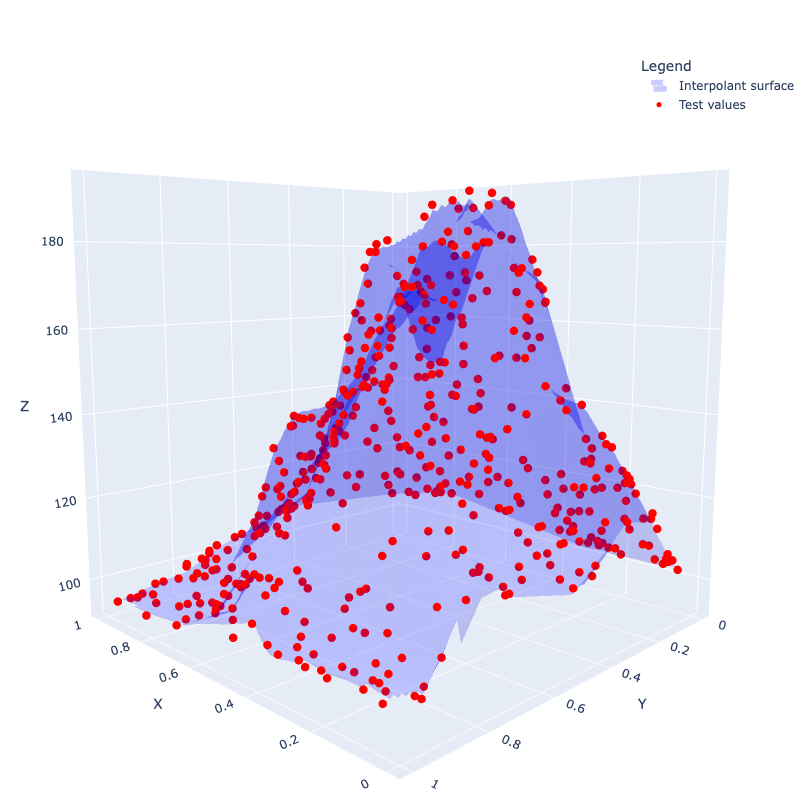}\hfill
\includegraphics[width=.49\textwidth]{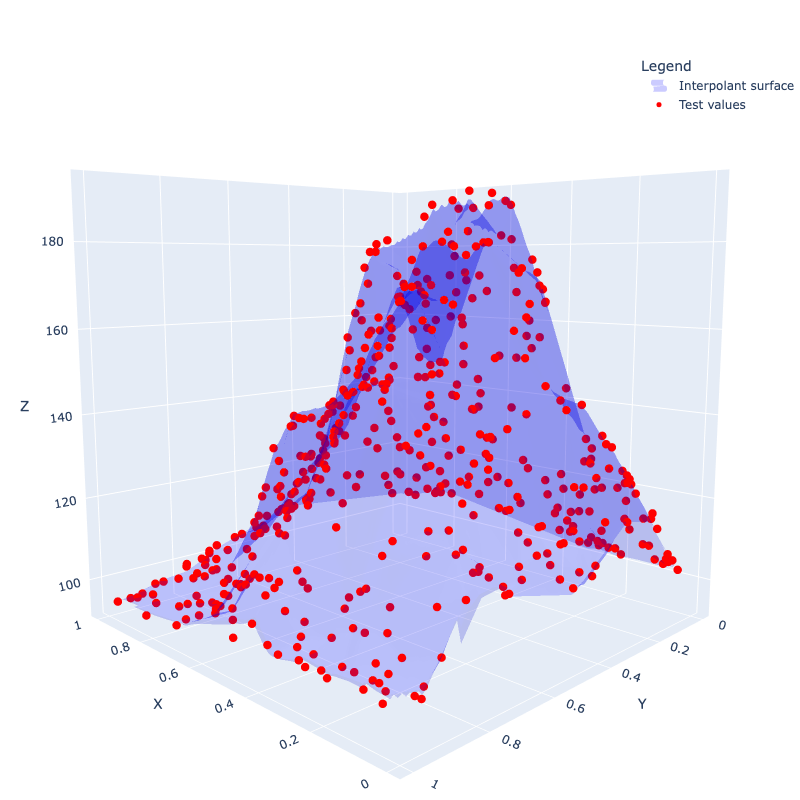}
    \caption{Interpolant surfaces and test values obtained using the M2 kernel (left) and the W2 kernel (right).}
    \label{fig:surfaces}
\end{figure}

To test the performance of the algorithm in the case of a real application, we use the \emph{volcano} dataset available in the statistical software package R \cite{R} representing 5307 elevation measurements obtained from Maunga Whau (Mt. Eden) in Auckland, NZ, sampled on a $10m \times 10m$ grid.
To produce an example with real data in line with the examples in the previous section, we consider a set of comparable dimensionality by selecting 
a random subset of 1500 points divided into a training set $(X, F)$ of $1000$ points and a test set  $(\Bar{X}, \Bar{F})$ of $500$ points (see Figure \ref{fig:train_test}). As explained in Section \ref{algorithms}, $(X, F)$ is further divided into $X_{train}, X_{val}$, $F_{train}, F_{val}$ when BO is used. 
For this experiment we focus on the interpolation case, i.e. setting $\Tilde{X} = X$ and so its divisions.  We consider a fixed parameter $\xi = 0.01$ and the M2 and W2 kernels. The results are shown in Table \ref{tab:real_data} and the obtained interpolant surfaces are displayed in Figure \ref{fig:surfaces}. We limited  ourselves to the use of the kernels M2 and W2 avoiding the GA  
because such data, coming from real measurements, have low regularity, while an essential requirement for the efficient use of the GA kernel is high regularity. Moreover, to make the results comparable we compute, alongside the $\mbox{MAE}_{\Bar{X}, \Bar{F}}(\cdot)$, its relative version, the  $\mbox{RMAE}_{\Bar{X}, \Bar{F}}(\cdot)$, obtained by dividing the  $\mbox{MAE}_{\Bar{X}, \Bar{F}}(\cdot)$ by the measurement with the highest absolute value.

\begin{table}[ht]
    \centering
    \begin{tabular}{|*{6}{p{14mm}|}}
        \hline
        kernel	& method	& time (s)	& $\mbox{MAE}_{\Bar{X}, \Bar{F}}$	& $\mbox{RMAE}_{\Bar{X}, \Bar{F}}$	& $\varepsilon^*$ \\ 
        \hline
        \multirow{2}{*}{M2} & LOOCV	& 9.66e+01 &	3.894942 &	0.020392	& 8.560000 \\
        \cline{2-6}
        \multirow{2}{*}{} & BO	& 6.38e+00	& 3.892905	& 0.020382 & 13.01788 \\
        \hline
        \multirow{2}{*}{W2} & LOOCV	& 1.04e+02 & 3.886987	& 0.020351	& 1.600000 \\
        \cline{2-6}
        \multirow{2}{*}{} & BO	& 8.24e+00 & 3.860059 & 0.020210 & 2.291239 \\
        \hline
    \end{tabular}
    \caption{Computational time, $\mbox{MAE}_{\Bar{X}, \Bar{F}}$, $\mbox{RMAE}_{\Bar{X}, \Bar{F}}$ and shape parameter using M2 and W2 kernels by LOOCV and BO.}
    \label{tab:real_data}
\end{table}

\section{Conclusions and Future Work}
\label{conclusion}

Comparing the parameter search in the case of interpolation with LOOCV and BO shows that the error values are almost always of the same magnitude, but what distinguishes the two methods is the calculation time, which in most cases is an order of magnitude lower with the BO approach. In the case of the approximation, we can see that there is a trade-off between computational cost and precision obtained. In particular, we can observe that, in general, considering a set of centers equal to $80\%$ of the number of nodes leads to lower computation effort without undermining the precision obtained to any great extent. Further clarifications on the results obtained must be made in the case in which a better error value is obtained with LOOCV. In these cases, the problem lies in the number of iterations of the BO; as the latter are limited, it can happen that the search gets bogged down in local minima without reaching the maximum. A larger number of iterations and an exploration-oriented parameter $\xi$ can solve the problem.  The performance of BO also depends on the choice of the covariance function, the prior distribution.  In the implementation given in \cite{BO}, the default covariance function is the Matérn $5/2$ kernel \cite{fas07}, which assumes that the function is quite smooth. A further insight could be changing the parameter of the covariance function or implementing different kernels for specific solution trying to improve performance even more.

 An additional extension of this work would be the application of Bayesian optimization to the Partition of Unity scheme \cite{Sandro} by locally determining the shape parameter for each subdomain.

\subsection*{Acknowledgments}  

The authors sincerely thank the reviewers for their constructive and valuable comments that enabled to significantly improve the paper. This research has been accomplished within the RITA \lq\lq Research ITalian network on Approximation\rq\rq\ and the UMI Group TAA \lq\lq Approximation Theory and Applications\rq\rq. This work has been supported by the INdAM--GNCS 2022 Project \lq\lq Computational methods for kernel-based approximation and its applications\rq\rq, code CUP$\_$E55F22000270001, by the 2022 Project \lq\lq Approximation
methods and models for life sciences\rq\rq\ and by the 2020 Project \lq\lq Mathematical methods in computational sciences\rq\rq\ funded by the Department of Mathematics \lq\lq Giuseppe Peano\rq\rq\ of the University of Torino. Moreover, the work of the first and second authors has been supported by the Spoke 1 ``Future HPC \& BigData'' of the Italian Research Center on High-Performance Computing, Big Data and Quantum Computing (ICSC) funded by MUR Missione 4 Componente 2 Investimento 1.4: Potenziamento strutture di ricerca e creazione di \lq\lq campioni nazionali di R$\&$S (M4C2-19)\rq\rq\ -- Next Generation EU (NGEU).  



\begin{thebibliography}{00}


\bibitem{review} L. Alzubaidi,  J. Zhang, A.J. Humaidi, Y. Duan, J. Santamaría, M.A. Fadhel, L. Farhan, Review of deep learning: Concepts, CNN architectures, challenges, applications, future directions, Journal of Big Data 8 (2021) 1--74. 

\bibitem{Biazar} B. Biazar, M. Hosami, An interval for the shape parameter in radial basis function approximation. Appl. Math. Comput. 315 (2017) 131--149.

\bibitem{Bishop} C.M. Bishop, Neural networks and their applications, Rev. Sci. Instrum. 65 (1994) 1803--1832.

\bibitem{Brochu} E. Brochu, V.M. Cora, N. De Freitas, A tutorial on Bayesian optimization of expensive cost functions, with application to active user modeling and hierarchical reinforcement learning, 2010, arXiv:1012.2599 

\bibitem{cav21a} R. Cavoretto, Adaptive radial basis function partition of unity interpolation: A bivariate algorithm for unstructured data, J. Sci. Comput. 87 (2021) 41.


\bibitem{cav19} R. Cavoretto, A. De Rossi, F. Dell’Accio, F. Di Tommaso, Fast computation of triangular Shepard interpolants, J. Comput. Appl. Math. 354 (2019) 457--470.

\bibitem{Sandro} R. Cavoretto, A. De Rossi, S. Lancellotti, E. Perracchione, Software implementation of the partition of unity method, Dolomites Res. Notes Approx. 15 (2022) 35--46. 

\bibitem{cav21} R. Cavoretto, A. De Rossi, M.S. Mukhametzhanov, Ya.D. Sergeyev, On the search of the shape parameter in radial basis functions using univariate global optimization methods, J. Global Optim. 79 (2021) 305--327.

\bibitem{cav22} R. Cavoretto, A. De Rossi, A. Sommariva, M. Vianello, RBFCUB: A numerical package for near-optimal meshless cubature on general polygons, Appl. Math. Lett. 125 (2022) 107704.

\bibitem{fas07} G.E. Fasshauer, Meshfree Approximation Methods with MATLAB, World Scientific, Singapore, 2007.

\bibitem{Fasshauer15}
G.E. Fasshauer, M.J. McCourt, Kernel-based Approximation Methods Using MATLAB, World Scientific, Singapore, 2015.

\bibitem{fas07b} G.E. Fasshauer, J.G. Zhang, On choosing \lq\lq optimal\rq\rq\ shape parameters for RBF approximation, Numer. Algorithms 45 (2007) 345--368.


\bibitem{Fornberg-Wright04} 
B. Fornberg, G. Wright, Stable computation of multiquadrics interpolants for all values of the shape parameter, Comput. Math. Appl. 47 (2004) 497--523.
 
\bibitem{Franke} R. Franke, Scattered data interpolation: tests of some methods, Math. Comp. 48 (1982) 181–-200.

\bibitem{Golberg} M.A. Golberg, C.S. Chen, S.R. Karur, Improved multiquadric approximation for partial differential equations, Eng. Anal. Bound. Elem. 18 (1997) 9--17.

\bibitem{Hardy} R.L. Hardy,  Multiquadric equations of topography and other irregular surfaces, J. Geophys. Res. 76 (1971) 1905--1915.

\bibitem{EI} D.R. Jones, M. Schonlau, W.J. Welch, Efficient Global Optimization of Expensive Black-Box Functions, J. Global Optim. 13 (1998) 455--492.


\bibitem{Larsson-Fornberg05} 
E. Larsson, B. Fornberg, Theoretical and computational aspects of multivariate interpolation with increasingly flat radial basis functions, Comput. Math. Appl. 49 (2005) 103--130.

\bibitem{lin22} L. Ling, F. Marchetti, A stochastic extended Rippa’s algorithm for LpOCV, Appl. Math. Letters 129 (2022) 107955.

\bibitem{Lizotte} D. Lizotte, Practical Bayesian Optimization, PhD thesis, University of Alberta, Edmonton, Alberta, Canada, 2008.

\bibitem{Marchetti} F. Marchetti, The extension of Rippa's algorithm beyond LOOCV, Appl. Math. Letters 120 (2021) 107262.

\bibitem{mir21} D. Mirzaei, The direct radial basis function partition of unity (D-RBF-PU) method for solving PDEs, SIAM J. Sci. Comput. 43 (2021) A54--A83. 


\bibitem{Mockus_1978} J. Mockus, V. Tiesis, A. Zilinskas, The application of Bayesian methods for seeking the extremum, Towards Global Optimization 2 (1978) 117--129.

\bibitem{BO} F. Nogueira, Bayesian optimization: Open source constrained global optimization tool for Python, \url{https://github.com/fmfn/BayesianOptimization}

\bibitem{sklearn} F. Pedregosa, G. Varoquaux, A. Gramfort, V. Michel, B. Thirion, O. Grisel, M. Blondel, P. Prettenhofer, R. Weiss, V. Dubourg, J. Vanderplas,A. Passos, D. Cournapeau, M. Brucher, M. Perrot, E. Duchesnay, Scikit-learn: Machine learning in Python, Journal of Machine Learning Research 12 (2011) 2825--2830.

\bibitem{R}
R Core Team (2020). R: A Language and Environment for Statistical Computing, R Foundation for Statistical Computing, Vienna, Austria.

\bibitem{Rasmussen} C.E. Rasmussen, C. Williams, Gaussian Processes for Machine Learning, MIT Press, 2006.

\bibitem{ren99} R.J. Renka, R. Brown, Algorithm 792: Accuracy test of ACM algorithms for interpolation of scattered data in the plane, ACM Trans. Math. Software 25
(1999) 78--94.

\bibitem{rip99} S. Rippa, An algorithm for selecting a good value for the parameter $c$ in radial basis function interpolation, Adv. Comput. Math. 11 (1999) 193--210.

\bibitem{rewiew_BO}
B. Shahriari, K. Swersky, Z. Wang, R.P. Adams,  N. de Freitas, Taking the human out of the loop: A review of Bayesian optimization,
Proceedings of the IEEE, 104-1 (2016) 148--175.

\bibitem{Practical} J. Snoek, H. Larochelle, R.P. Adams, Practical Bayesian Optimization of Machine Learning Algorithms, Advances in Neural Information Processing Systems 25 (2012) 2960--2968.

\bibitem{Scheuerer} M. Scheuerer, An alternative procedure for selecting a good value for the parameter c in RBF-interpolation, Adv. Comput. Math. 34 (2011)  105–126.

\bibitem{Trahan}
C.J. Trahan, R.W. Wyatt,  Radial basis function interpolation in the quantum trajectory method: optimization of the multi-quadric shape parameter, J. Comput. Phys. 185 (2003) 27--49. 

\bibitem{Uddin}
M. Uddin, On the selection of a good value of shape parameter in solving time-dependent partial differential equations using RBF approximation method, Appl. Math. Model. 38 (2014) 135--144.

\bibitem{scipy} P. Virtanen, R. Gommers, T.E. Oliphant, M. Haberland, T. Reddy, D. Cournapeau, E. Burovski, P. Peterson, W. Weckesser, J. Bright, S.J. van der Walt, M. Brett, J. Wilson, K. Jarrod Millman, N. Mayorov, A.R.J. Nelson, E. Jones, R. Kern, E. Larson, C.J. Carey, İ. Polat, Y. Feng, E.W. Moore, J. VanderPlas, D. Laxalde, J. Perktold, R. Cimrman, I. Henriksen, E.A. Quintero, C.R. Harris, A.M. Archibald, A.H. Ribeiro, F. Pedregosa, P. van Mulbregt, and SciPy 1.0 Contributors.  SciPy 1.0: Fundamental Algorithms for Scientific Computing in Python. Nature Methods 17(3) (2020)  261--272.

\bibitem{wen05}
H. Wendland, Scattered Data Approximation, Cambridge Monogr. Appl. Comput. Math., vol. 17, Cambridge Univ. Press, Cambridge, 2005.


\end{thebibliography}
\end{document}